\documentclass[twoside]{amsart} 
\usepackage{amssymb}
\usepackage{amsxtra}

\input xy
\xyoption{all}
\newdir{ >}{{}*!/-3mm/@{>}}
\newdir{ (}{{}*!/-3mm/@^{(}}

\theoremstyle{plain}
\newtheorem{proposition}{Proposition}[section]
\newtheorem{lemma}[proposition]{Lemma}
\newtheorem{corollary}[proposition]{Corollary}
\newtheorem{theorem}[proposition]{Theorem}
\newtheorem{property}[proposition]{Universal Property}

\theoremstyle{definition}
\newtheorem{notation}[proposition]{Notation}
\newtheorem{definition}[proposition]{Definition}
\newtheorem{remark}[proposition]{Remark}
\newtheorem{example}[proposition]{Example}

\allowdisplaybreaks

\DeclareMathOperator{\Ass}{Ass}
\DeclareMathOperator{\Spec}{Spec}
\DeclareMathOperator{\Gl}{Gl}
\DeclareMathOperator{\Hilb}{Hilb}
\DeclareMathOperator{\Hilbf}{\underline{Hilb}}
\newcommand{\Hf}{\underline{H}}
\DeclareMathOperator{\Char}{char}
\DeclareMathOperator{\satu}{sat} \newcommand{\sat}{{\satu}}
\DeclareMathOperator{\redu}{red} \newcommand{\red}{{\redu}}
\DeclareMathOperator{\Bin}{Bin}
\DeclareMathOperator{\init}{in}
\DeclareMathOperator{\Gin}{Gin}
\DeclareMathOperator{\lexu}{lex} \newcommand{\lex}{{\lexu}}

\newcommand{\sk}{\medskip}
\newcommand{\ra}{\rightarrow}
\newcommand{\Ra}{\Rightarrow}
\newcommand{\defiff}{:\Leftrightarrow}
\newcommand{\N}{\mathbb{N}}
\newcommand{\Z}{\mathbb{Z}}
\newcommand{\Q}{\mathbb{Q}}
\newcommand{\unipotent}{\mathcal{U}}
\newcommand{\Ideals}{\mathbb{I}}
\newcommand{\terms}[1][]{\mathbb{T}_{#1}}
\newcommand{\ideala}{{\mathfrak{a}}}
\newcommand{\idealb}{{\mathfrak{b}}}
\newcommand{\idealc}{{\mathfrak{c}}}
\newcommand{\ideald}{{\mathfrak{d}}}
\newcommand{\ideall}{{\mathfrak{l}}}
\newcommand{\idealp}{{\mathfrak{p}}}
\newcommand{\proj}[1][K]{{\mathbb{P}^n_{\!\!#1}}}
\newcommand{\sch}{{\mathcal{S}ch}}
\newcommand{\sets}{{\mathcal{S}et}}
\newcommand{\sheafo}{{\mathcal{O}}}
\newcommand{\sheafi}{{\mathcal{I}}}
\newcommand{\borelg}{>_\textup{\tiny{Bor}}}
\newcommand{\borelgeq}{\geq_\textup{\tiny{Bor}}}
\newcommand{\deglex}{\textup{hlex}}
\newcommand{\deglexg}{>_\textup{\tiny{hlex}}}
\newcommand{\degrevlex}{\textup{rlex}}
\newcommand{\degrevlexg}{>_\textup{\tiny{rlex}}}

\newcounter{numb}
\newenvironment{fact}%
{\begin{list}{}%
{\setlength{\labelwidth}{5.5ex} \setlength{\leftmargin}{6.5ex} \setlength{\labelsep}{1ex}
\setlength{\itemindent}{0ex} \setlength{\listparindent}{0ex}
\setlength{\parsep}{0ex} \setlength{\itemsep}{0ex}
\setlength{\topsep}{2ex} \setlength{\partopsep}{0ex}
\setlength{\parskip}{0ex}}
\addtocounter{numb}{1}\item[\textbf{(\thesection.\thenumb)}\hfill]\em}
{\end{list}}

\newenvironment{bulletlist}%
{\begin{list}{$\bullet$}%
{\setlength{\labelwidth}{2.0ex} \setlength{\leftmargin}{3.5ex} \setlength{\labelsep}{1.5ex}
\setlength{\itemindent}{0ex} \setlength{\listparindent}{0ex}
\setlength{\parsep}{0ex} \setlength{\itemsep}{0ex}
\setlength{\topsep}{1.0ex} \setlength{\partopsep}{0ex}
\setlength{\parskip}{0ex}}
}
{\end{list}}

\begin{document}

\title{Hilbert scheme strata defined by bounding cohomology}

\author{Stefan Fumasoli}
\date{Z\"urich, \today}

\subjclass[2000]{14C05, 14B15, 13P10}

\keywords{Hilbert scheme, local cohomology, Mall ideals, seqentially Cohen-Macaulayness}

\begin{abstract}
Let $\Hilb^p$ be the Hilbert scheme parametrizing the closed subschemes of $\proj$ with Hilbert polynomial $p\in\Q[t]$ over a field $K$ of characteristic zero. By bounding below the cohomological Hilbert functions of the points of $\Hilb^p$ we define locally closed subspaces of the Hilbert scheme. The aim of this paper is to show that some of these subspaces are connected. For this we exploit the edge ideals constructed by D.~Mall in \cite{M2000}. It turns out that these ideals are sequentially Cohen-Macaulay and that their initial ideals with respect to the reverse lexicographic term order are generic initial ideals.
\end{abstract}
\maketitle
\section{Introduction}
Let $\Hilb^p=\Hilb^p_{\proj}$ be the Hilbert scheme of projective space over a field $K$ with respect to a polynomial $p\in\Q[t]$. It is well known that $\Hilb^p$ is connected \cite{H1966}. For each point $x\in\Hilb^p$ let $h_x^0:\Z\ra\N$ denote the Hilbert function of the associated ideal sheaf $\sheafi^{(x)}\subset\smash{\sheafo_{\proj[\kappa(x)]}}$.
Several approaches have been chosen to prove that subsets of $\Hilb^p$ which are defined by bounding the functions $h_x^0$ are connected \cite{Go1988}, \cite{M2000}, \cite{P}, \cite{PS}.

For $i\in\N$ and $x\in\Hilb^p$ define the $i$th \emph{cohomological Hilbert function}
$$h^i_{x}:\Z\ra\N,\ m\mapsto\dim_{\kappa(x)}H^i\bigl(\proj[\kappa(x)],\sheafi^{(x)}(m)\bigr).$$
Fix a sequence $(f_i)_{i\in\N}$ of numerical functions $f_i:\Z\ra\N$. Then by the Semicontinuity Theorem 
$$H^{\geq}:=H^{\geq}_{\proj}:=\{x\in\Hilb^p\mid h^i_x\geq f_i\ \forall\,i\in\N\}$$ 
is a closed subspace of $\Hilb^p$ and 
$$H^{=}:=H^{=}_{\proj}:=\{x\in\Hilb^p\mid h^0_x=f_0,\ h^i_x\geq f_i\ \forall\,i\geq 1\}$$ 
is a locally closed subspace of $\Hilb^p$. In this paper we show (cf.\ Theorem \ref{connectedness by lines}):
  
\begin{fact}
$H^{\geq}$ and $H^{=}$ are connected if $\Char K=0$.
\end{fact}

We endow the subsets $H^{\geq}$ and $H^{=}$ with the induced reduced scheme structure. By flat base change and functorial arguments it follows that $H^{\geq}_{\proj[k]}\cong(H^{\geq}_{\proj}\times_K k)_\red$ and $H^{=}_{\proj[k]}\cong(H^{=}_{\proj}\times_K k)_\red$ for any field extension $K\subset k$. Hence, it suffices to show the connectedness of $H^{\geq}$ and $H^{=}$ in the case when $K$ is algebraically closed.

From now on we assume that $K$ is algebraically closed and of characteristic zero. The closed points of $\Hilb^p$ are precisely the saturated homogeneous ideals with Hilbert polynomial $q(t):=\binom{t+n}{n}-\nolinebreak p(t)$ of the polynomial ring $S:=K[X_0,\dots,X_n]$. For a homogeneous ideal $\ideala\subset S$ and $i\in\N$ define the $i$th \emph{locally cohomological Hilbert function} 
$$h^i_\ideala:\Z\ra\N,\ m\mapsto\dim_K H^i_{S_+}(\ideala)_m,$$ where $H^i_{S_+}(\ideala)$ denotes the (graded) $i$th local cohomology module of $\ideala$ with respect to the irrelevant ideal $S_+:=(X_0,\dots,X_n)$. Furthermore, let $h_\ideala$ denote the Hilbert function of $\ideala$. By the Serre-Grothendieck Correspondence the cohomology modules of sheaves correspond to local cohomology modules. 
Hence the set of closed points of $H^{=}$ equals
$$\Ideals^=:=\Big\{\ideala\subset S\Big|\ \parbox{6cm}{$\ideala$ is a saturated homogeneous ideal with $h_\ideala=f_0$ and $h^i_\ideala\geq f_{i-1}$ for all $i\geq 2$}\Big\}$$
and the set of closed points of $H^{\geq}$ equals
$$\Ideals^\geq:=\Big\{\ideala\subset S\Big|\ \parbox{8.3cm}{$\ideala$ is a saturated homogeneous ideal with Hilbert polynomial $q$ such that $h_\ideala\geq f_0$ and $h^i_\ideala\geq f_{i-1}$ for all $i\geq 2$}\Big\}.$$

A set $\Ideals$ of saturated homogeneous ideals of $S$ is said to be \emph{connected by Gr\"obner deformations} if for any two ideals $\ideala$, $\idealb\in\Ideals$ there exists a sequence of ideals $\ideala=$ $\idealc_1$, \dots, $\idealc_r=\idealb$ in $\Ideals$ such that $\idealc_i$ is the saturation of the initial ideal or of the the generic initial ideal of $\idealc_{i+1}$ with respect to some term order or vice versa for all $i\in\{1,\dots,r-1\}$. We exploit the ideals constructed by D.~Mall in \cite{M2000} to show that $\Ideals^=$ and $\Ideals^\geq$ are connected by Gr\"obner deformations. With techniques, described in Chapter 15.8 of \cite{E}, one shows that connectedness by Gr\"obner deformations implies connectedness in sense of topology for the sets $H^{\geq}$ and $H^{=}$ (cf.\ proof of Theorem~\ref{connectedness by lines}).

To prove that a set $\Ideals$ of saturated homogeneous ideals of $S$, defined by bounding the locally cohomological Hilbert functions, is connected by Gr\"obner deformations, \emph{sequentially Cohen-Macaulay} ideals (i.e. ideals $\ideala\subset S$ such that $S/\ideala$ is a sequentially Cohen-Macaulay $S$-module) play an important role: 

\begin{fact}
Let $\ideala\subset S$ be a homogeneous ideal, and let $\Gin_\degrevlex\ideala$ be its generic initial ideal with respect to the homogeneous reverse lexicographic order. Then $h^i_\ideala=h^i_{\Gin_\degrevlex\ideala}$ for all $i\in\N$ if and only if $\ideala$ is sequentially Cohen-Macaulay (cf.\ \cite{HS}).
\end{fact}
\noindent
On the other hand we always have:
\begin{fact}
Let $\ideala\subset S$ be a homogeneous ideal. Then $h^i_\ideala\leq h^i_{\init\ideala}$ for all $i\in\N$ with respect to any term order on $S$ (cf.\ \cite{S}).
\end{fact}

In his Habilitations\-schrift of 1997, D. Mall showed that the two sets 
$$\Ideals_=:=\{\ideala\subset S\mid\text{$\ideala$ is a saturated homogeneous ideal with $h_\ideala=f_0$}\},$$  
$$\Ideals_{\geq}:=\{\ideala\subset S\mid\text{homogeneous, saturated, with Hilb.\ poly.\ $q$ and $h_\ideala\geq f_0$}\}$$ 
are connected by Gr\"obner deformations. Observe that $\Ideals^=\subset\Ideals_=$, $\Ideals^\geq\subset\Ideals_\geq$ are subsets. Indeed Mall described an algorithm providing the following fact: 
\begin{fact}\label{pm}
Assume that\/ $\Ideals_=$ is not empty. Then there exists an ideal $\ideall^{f_0}\in\Ideals_=$, (a so called \emph{growth-height-lexicographic ideal} depending uniquely on $f_0$, cf.\ \cite{M1997}), such that for each $\ideala\in\Ideals_=$ there exists a sequence of ideals $\idealc_1$, \dots, $\idealc_r$ in $\Ideals_=$ such that 
\begin{bulletlist}
\item $\init_\degrevlex\idealc_1=\Gin_\degrevlex\ideala$, 
\item $\init_\degrevlex\idealc_i=\init_\deglex\idealc_{i-1}$ for all $i\in\{2,\dots,r\}$,
\item $\init_\deglex\idealc_r=\ideall^{f_0}$.
\end{bulletlist}
Moreover, there exists a sequence of ideals $\ideald_1$, \dots, $\ideald_s$ in $\Ideals_\geq$ such that
\begin{bulletlist}
\item $\init_\degrevlex\ideald_1=\ideall^{f_0}$, 
\item $\init_\degrevlex\ideald_i=(\init_\deglex\ideald_{i-1})^\sat$ for all $i\in\{2,\dots,s\}$, 
\item $(\init_\deglex\ideald_s)^\sat=\ideall^p$, 
\item $h_{\init_\degrevlex\ideald_i}=h_{\ideald_i}\leq h_{(\init_\deglex\ideald_i)^\sat}$ for all $i\in\{1,\dots,s\}$.
\end{bulletlist}
where $\ideall^q$ denotes the unique saturated lexicographic ideal with Hilbert polynomial $q$ (cf.\ \cite{M2000}).
\end{fact}
I will call the ideals $\idealc_i$ and $\ideald_i$ \emph{Mall ideals}. These ideals have a lot of nice properties: At first, they are generated by monomials and by homogeneous binomials. More precisely, the binomials are all parallel. Hence Mall ideals are \emph{edge ideals} (cf.\ \cite{AS}): For every Mall ideal $\idealc$ there exists $\rho\in\Z^{n+1}\setminus\{0\}$ with $\rho_0+\dots+\rho_n=0$ such that $\idealc$ is homogeneous with respect to the induced $\Z^{n+1}/\rho\Z$-grading of $S$. In particular they are \emph{edge providing}: Any Mall ideal has exactly two initial ideals, namely one with respect to the homogeneous lexicographic order $\geq_\deglex$ and one with respect to the homogeneous reverse lexicographic order $\geq_\degrevlex$. Moreover, these two initial ideals are \emph{fixed by the action of the Borel group}.

In the main part of this paper we prove the following two statements (cf.\ Theorem \ref{Mall ideals are sCM} and Theorem \ref{Gin = in for Mall ideals}):

\begin{fact}
Mall ideals are sequentially Cohen-Macaulay.
\end{fact}

\begin{fact}
If $\idealc\subset S$ is a Mall ideal, then $\Gin_\degrevlex\idealc=\init_\degrevlex\idealc$.
\end{fact}
\noindent
Therefore, on use of the facts (0.2), (0.3) and (0.4) it follows (cf.\ Proposition \ref{connectedness by Groebner def}):
\begin{fact}
The sets $\,\Ideals^=$ and $\,\Ideals^\geq$ are connected by Gr\"obner deformations.
\end{fact}
\noindent
This finishes the proof of statement (0.1).

The calculations of the examples (cf.~\ref{example1}, \ref{example2}, \ref{counterexample}) have been performed with CoCoA \cite{Cocoa}.

This paper is based on my Ph.~D.~thesis \cite{F} to which I refer for technical details.

\sk
\section{Preliminaries}
In this section we introduce the needed combinatorial and algebraic notions and collect some results about initial and generic initial ideals. The latter ones play an important role in the proof of connectedness of Hilbert schemes: A homogeneous ideal $\ideala$ in a polynomial ring $K[X_1,\dots,X_n]$ over a field $K$ and its generic initial ideal $\Gin\ideala$ have the same Hilbert function. Furthermore, by means of weight orders, they are connected in the Hilbert scheme by a sequence of lines (cf.~section~\ref{Hilbert scheme}).

Generic initial ideals are Borel-fixed, which means that they are fixed under the action of upper triangular matrices. If $K$ has characteristic zero, they correspond to Borel sets which have a combinatorial behavior which is rather easy to understand (cf.~section~\ref{subsection Borel sets}). 
\subsection{Notations and definitions}\label{subsection notations}
\begin{notation}
Let $\N$ denote the set of nonnegative integers. Throughout this paper $n$, $d\in\N\setminus\{0\}$ are two positive integers.

Let $K$ be a field with $\Char K=0$. 

We fix a polynomial ring $S:=K[X_1,\dots,X_n]$. (In the last section $S$ will denote the polynomial ring $K[X_0,\dots,X_n]$ in one variable more.) Any polynomial ring is endowed with the standard $\Z$-grading. 

Let denote $\terms$ the set of monomials of $S$. Define $\N_d^n:=\{a\in\N^n\mid a_1+\dots+a_n=d\}$. Then there is a bijection $\N^n_d\ra\terms\cap S_d,\ a\mapsto X^a:=X_1^{a_1}\cdots X_n^{a_n}.$ Via this bijection we may identify the elements of $\N^n$ with the monomials in $\terms$.

Let denote $\deglexg$ the \emph{homogeneous lexicographic order} and $\degrevlexg$ the \emph{reverse lexicographic order} of $\terms$. Both term orders are \emph{admissible}, i.e. $X_1>\dots>X_n$ and $\deg u>\deg v\ \Rightarrow\ u>v$ for all $u,v\in\terms$.
\end{notation}

\begin{definition}
Let $B\subset\N^n_d$. A subset $A\subset B$ is called a \emph{lexicographic segment} of $B$ if for all $a,b\in B$ with $a\in A$ and $b\deglexg a$ we have $b\in A$.

The subset $B\subset\N^n_d$ is called \emph{lexicographic} if $B$ is a lexicographic segment of $\N^n_d$. In this case, also the set $X^B:=\{X^b\in\terms\mid b\in B\}$ is called \emph{lexicographic}.

A monomial ideal $\ideala\subset S$ is called a \emph{lex ideal} if $\ideala_i\cap\terms$ is lexicographic for all $i\in\N$.
\end{definition}

\begin{definition}
Let $R=\bigoplus_{i\in\N}R_i$ be a homogeneous Noetherian ring and $M$ a finitely generated graded $R$-module. Let $R_+:=\bigoplus_{i>0}R_i$ denote the \emph{irrelevant ideal} of $R$. For $i\in\N$ let $H^i_{R_+}(M)$ denote the $i$th \emph{local cohomology module} of $M$ with respect to $R_+$, endowed with its natural grading (cf.\ \cite[Chap.~12]{BS}).

Let denote by $h_M:\Z\ra\N,\ i\mapsto\dim_K M_i$ the \emph{Hilbert function} of $M$. 

A polynomial $p\in\Q[t]$ is called an \emph{admissible Hilbert polynomial} if there exists a homogeneous ideal $\ideala\subset S$ with Hilbert polynomial $p$.
\end{definition}

\begin{remark}\label{existence of lex ideals}
Let $\ideala\subset S$ be a homogeneous ideal. Then there exists a unique lex ideal $\ideala^\lex\subset S$ with $h_\ideala=h_{\ideala^\lex}$ (cf.~\cite[4]{Sp}).

Let $p\in\Q[t]$ be an admissible Hilbert polynomial. Then there exists a unique saturated lex ideal $\ideall^p\subset S$ with Hilbert polynomial $p$ (cf.~\cite[1.11]{F}).
\end{remark}

\begin{remark}\label{ideal and initial ideal have same Hilbert function}
Let $\ideala\subset S$ be a homogeneous ideal and let $\tau$ be a term order of\/ $\terms$. Then $h_\ideala=h_{\init_\tau\ideala}$ (cf.~\cite[15.26]{E}).
\end{remark}

We introduce some more combinatorial notations:

\begin{notation}\label{combinatorial notations}
For any set $L$ and any integer $m\in\N$ let $L^{(n,m)}$ denote the set of all $n\times m$ matrices $[M_{ij}\mid 1\leq i\leq n,\:1\leq j\leq m]$ with entries $M_{ij}\in L$. 

Set $U(n):=\{M\in \N^{(n,n)}\mid M_{ij}=0\ \ \forall\ 1\leq j<i\leq n\}$. 

Let $g=[g_{ij}\mid 1\leq i\leq n,\:1\leq j\leq n]\in K^{(n,n)}$ be a matrix. Then we denote by $g:S\ra S$ the homomorphism of $K$-algebras defined by $X_j\mapsto\sum_{i=1}^n g_{ij}X_i$ for $1\leq j\leq n$.

\sk
For two functions $f,g\colon\Z\ra\N$ we write $f\geq g$ if $f(k)\geq g(k)$ for all $k\in \Z$, and we write $f>g$ if $f\geq g$ and if there exists $k\in\Z$ such that $f(k)>g(k)$.

For $1\leq i\leq n$, let $e_i\in\N^n$ denote the standard vector with $(e_i)_j=
1$ if $i=j$ and $(e_i)_j=0$ otherwise.

\sk
For the following notations let $a=(a_1,\dots,a_n)\in\Z^n$ and $A\subset\N^n$:

$m(a):=\max{(\{1\leq i\leq n\mid a_i\neq 0\}\cup\{1\})}$,

$A+a:=\{b+a\mid b\in A\}$,

$a^*:=a-a_n e_n=(a_1,\dots,a_{n-1},0)$,

$A^*:=\{b^*\mid b\in A\}$,

$a^+:=(\max{\{a_1,0\}},\dots,\max{\{a_n,0\}})\in\N^n$, 

$a^-:=(-\min{\{a_1,0\}},\dots,-\min{\{a_n,0\}})\in\N^n$, so that $a=a^+-a^-$. 
\end{notation}

\subsection{Borel sets}\label{subsection Borel sets}
Since the characteristic of $K$ is zero, the Borel-fixed ideals $\ideala\subset S$ are monomial ideals which are characterized by the following property: If a monomial $m\in\ideala$ is divisible by an indeterminate $X_j$, then $\frac{X_i}{X_j}\,m\in\ideala$ for all $i\leq j$. In each homogeneous component they correspond to so called Borel sets. Borel sets are the Borel order analogue to lexicographic sets. 

There are several equivalent ways to define the Borel order. The most plausible is the following one: For all monomials $m\in\terms$ and for all $1\leq k<n$ set $X_k\,m\borelg X_{k+1}\,m$ and take the associative hull:
\begin{definition}\label{def - Borel order, set - def}
Define the \emph{Borel order} $\borelgeq$ of $\N^n_d$ as follows. Let $a,b\in\N_d^n$. 

$a\borelgeq b\ \defiff\ \exists\,\alpha_1,\dots,\alpha_{n-1}\in\N:a-b=\sum_{j=1}^{n-1}\alpha_j(e_j-e_{j+1})$.

A set $A\subset\N^n_d$ is called a \emph{Borel set} if for all $a,b\in\N^n_d$ with $a\in A$ and $b\borelgeq a$ we have $b\in A$. In this case, also the set $X^A\subset\terms$ is called a \emph{Borel set}.
\end{definition}

For technical reasons it will be more convenient to have a description of the Borel order by upper triangular integer matrices (s.\ Lemma \ref{properties of Borel order}).

As a consequence of the fact that generic initial ideals are Borel-fixed independently of the admissible term order (Proposition \ref{existence of gin}) we have the following characterization of the Borel order:

\begin{lemma}[{\cite[2.2]{C}}]\label{Borel greater implies greater for any term order}
Let $a$, $b\in\N^n_d$. Then it holds $a\borelgeq b$ if and only if $X^a\geq X^b$ for all admissible term orders $\geq$ of\/ $\terms$.
\end{lemma}

\begin{remark}\label{linearity of Borel order}
Let $a$, $b\in\N^n_d$ and $\rho\in\Z^n$ such that $a+\rho$, $b+\rho\in\N^n_d$. Then $a\borelgeq b$ if and only if $a+\rho\borelgeq b+\rho$.
\end{remark}

\begin{definition}\label{def - monomial, Borel ideal - def}
A monomial ideal $\ideala\subset S$ is called a \emph{Borel ideal} if $\ideala_i\cap\terms$ is a Borel set for all $i\in\N$.
\end{definition}

\begin{remark}\label{Groebner bases of Borel ideals}
If $\ideala\subset S$ is a Borel ideal, then  $\ideala^\sat=(\ideala:X_n^\infty)$ (cf.~\cite[15.24]{E}). Hence, if $B\subset\N^n_d$ is a Borel set, then the set $X^{B^*}$ is a Gr\"obner basis of $(X^B)^\sat$ with respect to any term order of\/ $\terms$.
\end{remark}

\begin{lemma}\label{properties of Borel order}
Let $a,b\in\N_d^n$. Then the following are equivalent:
\begin{align*}
&\textup{(i)}&&a\borelgeq b,\\
&\textup{(ii)}&&\exists M\in U(n):\ \sum_{i=1}^n M_{ij}=b_j\ \ \forall\ 1\leq j\leq n,\ \sum_{j=1}^n M_{ij}=a_i\ \ \forall\ 1\leq i\leq n.
\end{align*}
\end{lemma}

\begin{proof}
(i)$\Ra$(ii): Let $a\borelgeq b$. Choose $\alpha_k\in\N$ for all $1\leq k<n$  such that $a-b=\sum_{k=1}^{n-1}\alpha_k(e_k-e_{k+1})$. Set $\alpha_0:=0$ and 
$m:=-\sum_{i=1}^n \min{\{0,b_i-\alpha_{i-1}\}}.$ We will construct a sequence of matrices $M(0)$, \dots, $M(m)\in\Z^{(n,n)}$ such that for all $0\leq k\leq m$ the following properties hold:
\begin{align*}
&\text{(1)}\hspace{12pt}M(k)_{ij}=0\text{ for all }1\leq j<i\leq n,\\
&\text{(2)}\hspace{12pt}M(k)_{ij}\geq 0\text{ for all }1\leq i<j\leq n,\\
&\text{(3)}\hspace{12pt}\sum_{i=1}^n M(k)_{ij}=b_j\text{ for all }1\leq j\leq n,\\
&\text{(4)}\hspace{12pt}\sum_{j=1}^n M(k)_{ij}=a_i\text{ for all }1\leq i\leq n,\\
&\text{(5)}\hspace{12pt}\sum_{i=1}^n \min{\{0,M(k)_{ii}\}}=k-m.
\end{align*}
Property (5) implies that $M(m)_{ii}\geq 0$ for all $1\leq i \leq n$, hence $M(m)$ is the requested matrix in $U(n)$.

Define $M(0)\in\Z^{(n,n)}$ by $M(0)_{ij}:=
\begin{cases}
\alpha_{i},&\text{if $i=j-1$;}\\
b_i-\alpha_{i-1},&\text{if $i=j$;}\\
0,&\text{otherwise.}
\end{cases}$\\
It is clear that $M(0)$ has the properties (1), (2), (3) and (5). Setting $\alpha_n:=0$, we have for all $1\leq i\leq n$
$$\sum_{j=1}^n M(0)_{ij}=b_i-\alpha_{i-1}+\alpha_i=b_i+\bigl(\sum_{k=1}^{n-1}\alpha_k(e_k-e_{k+1})\bigl)_i=a_i,$$
whence $M(0)$ has also property (4).

If $m>0$, we construct $M(1)$, \dots, $M(m)$ recursively. Let $0\leq k<m$ and assume that $M(k)\in\Z^{(n,n)}$ with the required properties is constructed already. Since $k<m$, by property (5) there exists 
$$l:=\min{\{1\leq j\leq n\mid M(k)_{jj}<0\}}.$$ 
Properties (1), (3) and (4) imply that there exist $p,q\in\N$ with $1\leq p<l<q\leq n$ and $M(k)_{pl}$, $M(k)_{lq}>0$. Now define $M(k+1)\in\Z^{(n,n)}$ by $$M(k+1)_{ij}:=
\begin{cases}
M(k)_{ij}+1,&\text{if $(i,j)\in\{(l,l),(p,q)\}$;}\\
M(k)_{ij}-1,&\text{if $(i,j)\in\{(p,l),(l,q)\}$;}\\
M(k)_{ij},&\text{otherwise.}
\end{cases}$$
It is clear that $M(k+1)$ has the required properties.

\sk
(ii)$\Ra$(i): For $1\leq k<n$ set $\alpha_k:=\sum_{i=1}^k a_i-b_i$. If there exists $M\in U(n)$ such that $\sum_{i=1}^n M_{ij}=b_j$ for all $1\leq j\leq n$ and $\sum_{j=1}^n M_{ij}=a_i$ for all $1\leq i\leq n$, then it follows that $\alpha_k\in\N$ for all $1\leq k<n$. Furthermore, an easy computation gives
$a-b=\sum_{k=1}^{n-1}\alpha_k(e_k-e_{k+1})$ (cf.~\cite[1.26]{F}).
\end{proof}

\subsection{Generic initial ideals and reverse lexicographic order}\label{subsection Gin}
Generic initial ideals have a lot of nice properties. We already mentioned that they are fixed under the action of the Borel group. In this section we show that if $S$ is endowed with the reverse lexicographic order, then their formation commutes with saturation. (As before, we assume that $\Char K=0$). 

\begin{definition}
The \emph{unipotent subgroup} $\unipotent\subset\Gl(n,K)$ is the group of all upper triangular matrices with ones on the diagonal.
\end{definition}

\begin{proposition}[{\cite[Chap. 15.9]{E}}]\label{existence of gin}
Let $\ideala\subset S$ be a homogeneous ideal and $\tau$ an admissible term order of\/ $\terms$. 

a) There is a non-empty Zariski open set $U\subset\Gl(n,K)$ and a unique ideal\/ $\Gin_\tau\ideala\subset S$ such that $\Gin_\tau\ideala=\init_\tau g(\ideala)$ for all $g\in U$. Furthermore, the open set $U$ meets the unipotent group $\unipotent$.

The ideal\/ $\Gin_\tau\ideala$ is called the \emph{generic initial ideal} of $\ideala$ with respect to $\tau$.

b) The generic initial ideal $\Gin_\tau\ideala$ is \emph{Borel-fixed}, i.\,e.\ for all upper triangular matrices $g\in\Gl(n,K)$ it holds $g(\Gin_\tau\ideala)=\Gin_\tau\ideala$. 

c) An ideal $\idealb\subset S$ is Borel-fixed if and only if it is a Borel ideal.
\end{proposition}

\begin{corollary}\label{Gin=in for unipotent fixed ideals}
If a homogeneous ideal of $S$ remains fixed under the action of $\unipotent$, then its generic initial ideal and its initial ideal with respect to any admissible term order coincide.
\end{corollary}

\begin{proposition}\label{Gin and sat commute}
Let $\ideala\subset S$ be a homogeneous ideal. Then
$$\Gin_\degrevlex(\ideala^\sat)=(\Gin_\degrevlex\ideala)^\sat.$$
\end{proposition}

\begin{proof}
Let $P:=\bigcup_{\idealp\in\Ass(S/\ideala)\setminus\{S_+\}}\idealp$ be the set of all elements of $S$ contained in some relevant associated prime of $S/\ideala$. We first show that $(\ideala: u^\infty)=\ideala^\sat$ for all $u\in S_1\setminus P$. Let $u\in S_1\setminus P$. Let $r\in\N$ be such that $\ideala^\sat=(\ideala: S_+^r)$. Since $u^r\notin P$, it holds $(0:_{S/\ideala}u^r)=H^0_{S_+}(0:_{S/\ideala}u^r)$ (cf.\ \cite[18.3.8\,(iii)]{BS}). It follows $(\ideala: u^r)=\ideala^\sat$.

We next prove that there exists a non-empty Zariski open set $U\subset\Gl(n,K)$ such that $g(\ideala^\sat)=(g(\ideala):X_n^\infty)$ for all $g\in U$. Since $\#K=\infty$, the open subset $S_1\setminus P\subset S_1$ is not empty, whence $U:=\{g\in\Gl(n,K)\mid g^{-1}(X_n)\notin P\}$ 
is a dense open subset of $\Gl(n,K)$. Let $g\in U$. Then it holds $g(\ideala^\sat)=g(\ideala:\nolinebreak g^{-1}(X_n)^\infty)=(g(\ideala):X_n^\infty)$.

Choose $g\in U$ such that $\Gin_\degrevlex\ideala=\init_\degrevlex g(\ideala)$ and $\Gin_\degrevlex(\ideala^\sat)=\init_\degrevlex g(\ideala^\sat)$. Since $\Gin_\degrevlex(\ideala)$ is a Borel ideal, $(\Gin_\degrevlex \ideala:X_n^\infty)=(\Gin_\degrevlex\ideala)^\sat$ (cf.\ Remark~\ref{Groebner bases of Borel ideals}). By \cite[15.12]{E} we have $\init_\degrevlex(g(\ideala):\nolinebreak X_n^\infty)=(\init_\degrevlex g(\ideala):\nolinebreak X_n^\infty)$. Altogether we obtain
$\Gin_\degrevlex(\ideala^\sat)=\init_\degrevlex g(\ideala^\sat)=\init_\degrevlex(g(\ideala):\nolinebreak X_n^\infty)=(\init_\degrevlex g(\ideala):\nolinebreak X_n^\infty)=(\Gin_\degrevlex \ideala:\nolinebreak X_n^\infty)=(\Gin_\degrevlex\ideala)^\sat.$
\end{proof}

\sk
\section{Mall ideals}\label{section binomial ideals}
Fix a polynomial $q\in\Q[t]$ and a function $f:\Z\ra\N$. Define the sets
$$\Ideals_=:=\{\ideala\subset S\mid\text{$\ideala$ is a saturated homogeneous ideal with $h_\ideala=f$}\},$$  
$$\Ideals_{\geq}:=\Big\{\ideala\subset S\ \Big|\ \parbox{7.2cm}{$\ideala$ is a saturated homogeneous ideal with Hilbert polynomial $q$ and with $h_\ideala\geq f$}\Big\}.$$ 
D.~Mall described an algorithm which yields the following proposition (cf.\ the picture below): 

\begin{proposition}[\cite{M2000}]\label{Mall's main theorem}
Assume that $\Ideals_=$ is not empty. Then there exists an ideal $\ideall^f\in\Ideals_=$, (a so called \emph{growth-height-lexicographic ideal} depending uniquely on~$f$, cf.~\cite{M1997}), such that for each $\ideala\in\Ideals_=$ there exists a sequence of ideals $\idealc_1$,~\dots,~$\idealc_r$ in $\Ideals_=$ such that 
\begin{bulletlist}
\item $\init_\degrevlex\idealc_1=\Gin_\degrevlex\ideala$, 
\item $\init_\degrevlex\idealc_i=\init_\deglex\idealc_{i-1}$ for all $i\in\{2,\dots,r\}$,
\item $\init_\deglex\idealc_r=\ideall^f$.
\end{bulletlist}

Moreover, there exists a sequence of ideals $\ideald_1$, \dots, $\ideald_s$ in $\Ideals_\geq$ such that
\begin{bulletlist}
\item $\init_\degrevlex\ideald_1=\ideall^f$, 
\item $\init_\degrevlex\ideald_i=(\init_\deglex\ideald_{i-1})^\sat$ for all $i\in\{2,\dots,s\}$, 
\item $(\init_\deglex\ideald_s)^\sat=\ideall^q$, 
\item $h_{\init_\degrevlex\ideald_i}=h_{\ideald_i}<h_{(\init_\deglex\ideald_i)^\sat}$ for all $i\in\{1,\dots,s\}$.
\end{bulletlist}
\end{proposition}

\begin{center}
\framebox{$\def\objectstyle{\scriptstyle}\def\labelstyle{\scriptstyle}\xymatrix@C=-2pt@R=7pt{
*i{\txt{100\\2}}&&&&\txt{\ \\$\Ideals_=$}&&&&&&\txt{\ \\$\Ideals_\geq$}\\
\ideala\ar@{-}[dr]
&&\idealc_1\ar@{-}[dl]\ar@{-}[dr]
&&\idealc_2\ar@{-}[dl]\ar@{.}[dr]
&\dots
&\idealc_r\ar@{.}[dl]\ar@{-}[dr]
&&\ideald_1\ar@{-}[dl]\ar@{-}[dr]
&&\ideald_2\ar@{-}[dl]\ar@{.}[dr]
&\dots
&\ideald_s\ar@{.}[dl]\ar@{-}[dr]\\
&\txt{$\scriptstyle\Gin_\degrevlex\ideala$\\$\scriptstyle=\init_\degrevlex\idealc_1$}
&&\txt{$\scriptstyle{\init_\deglex\idealc_1}$\\$\scriptstyle=\init_\degrevlex\idealc_2$}
&&*+<8ex,5.2ex>{}
&&\txt{$\scriptstyle\init_\deglex\idealc_r=$\\\hspace*{-7pt}$\scriptstyle\ideall^f=\init_\degrevlex\ideald_1$\hspace*{-7pt}}
&*+<2ex,5.2ex>{}
&\txt{\hspace*{-6pt}$\scriptstyle{(\init_\deglex\ideald_1)^\sat}$\hspace*{-14pt}\\$\scriptstyle=\init_\degrevlex\ideald_2$}
&&*+<8ex,5.2ex>{}
&&\txt{\hspace*{-10pt}$\scriptstyle(\init_\deglex\ideald_s)^\sat$\hspace*{-10pt}\\$\scriptstyle=\ideall^q$}
\save "1,1"."3,9"*[F.]\frm{}\restore
}$\hspace{10pt}}
\end{center}
\medskip

Hence the two sets $\Ideals_=$ and $\Ideals_\geq$ are connected by Gr\"obner deformations in the following sense:

\begin{definition}
A set $\Ideals$ of saturated homogeneous ideals of $S$ is said to be \emph{connected by Gr\"obner deformations} if for any two ideals $\ideala$, $\idealb\in\Ideals$ there exists a sequence of ideals $\ideala=\idealc_1$, \dots, $\idealc_r=\idealb$ in $\Ideals$ such that for all $i\in\{1,\dots,r-1\}$ it holds $\idealc_i=(\init_\tau\idealc_{i+1})^\sat$, $\idealc_i=(\Gin_\tau\idealc_{i+1})^\sat$, $\idealc_{i+1}=(\init_\tau\idealc_i)^\sat$, or $\idealc_{i+1}=(\Gin_\tau\idealc_i)^\sat$ with respect to some term order $\tau$ on $S$. 
\end{definition}

\begin{definition}
An ideal $\idealc\subset S$ is called \emph{Mall ideal}, if there exist $f:\Z\ra\N$, $q\in Q[t]$, $\ideala\in\Ideals_=$, $r$, $s\in\N$ and $\idealc_1$,~\dots,~$\idealc_r$, $\ideald_1$, \dots, $\ideald_s$, computed by Mall's algorithm as in Proposition \ref{Mall's main theorem}, such that $\idealc\in\{\idealc_1,\dots,\idealc_r,\ideald_1,\dots,\ideald_s\}$.
\end{definition}

\begin{example}\label{example1}
Let $S:=K[x,y,z,t]$ and $q(t):=\binom{t+3}{3}-(2t+3)=\frac{t^3}{6}+t^2-\frac{t}{6}-2$. Define $f:\Z\ra\N$ by $f(m):=0$ for $m\leq 2$ and $f(m):=q(m)$ for $m\geq 3$. The growth-height-lexicographic ideal with respect to $f$ is $\ideall^f=(x^2,xy,xz,y^3,y^2z)\in\Ideals_=.$ The saturated lex ideal with Hilbert polynomial $q$ is $\ideall^q=(x,y^3,y^2z^2)\in\Ideals_\geq.$

There are just three saturated Borel ideals with Hilbert polynomial $q$, namely $\ideall^f$, $\ideall^q$ and $\idealb:=(x^2,xy,y^2,xz^2)\in\Ideals_=.$

There exists a Mall ideal $\idealc=(y^2-xz,x^2,xy,xz^2)\in\Ideals_=$ such that $\init_\degrevlex\idealc=\idealb$ and $\init_\deglex\idealc=\ideall^f$. Moreover, there exists a Mall ideal $\ideald=(y^2z-xt^2,x^2,xy,xz,y^3)\in\Ideals_=$ such that $\init_\degrevlex\ideald=\ideall^f$ and $(\init_\deglex\ideald)^\sat=\ideall^q$. Hence, we get the following picture:
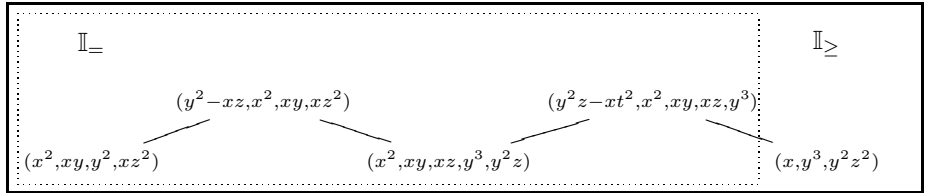
\begin{center}
\framebox{$$\def\objectstyle{\scriptstyle}\def\labelstyle{\scriptstyle}\xymatrix@C=0pt@R=7pt{
\txt{\hspace*{50pt}\\$\Ideals_=$}&&&&\txt{\ \\$\Ideals_\geq$}\\
&(y^2-xz,x^2,xy,xz^2)\ar@{-}[dl]\ar@{-}[dr]
&&(y^2z-xt^2,x^2,xy,xz,y^3)\ar@{-}[dl]\ar@{-}[dr]\\
(x^2,xy,y^2,xz^2)
&&(x^2,xy,xz,y^3,y^2z)
&*+<18.5ex,4.0ex>{}
&(x,y^3,y^2z^2)
\save "1,1"."3,4"*[F.]\frm{}\restore
}$$\hspace{10pt}}
\end{center}
\end{example}

Mall ideals have a lot of nice properties. In particular, they are generated by monomials and by binomials. The generating binomials are homogeneous with respect to the standart $\Z$-grading of $S$, and furthermore they are parallel. That means that there exists $\rho\in\Z^n$ such that they are homogeneous with respect to the $\Z^n/\rho\Z$-grading of $S$. Hence Mall ideals are edge ideals in the sense of \cite{AS}. In section \ref{subsection seq. CM} we show that Mall ideals are sequentially Cohen-Macaulay. In section \ref{subsection generic initial ideal} we show that the initial ideal of a Mall ideal coincides with its generic initial ideal. For those two facts we need some combinatorial properties of Mall ideals which are stated in section \ref{subsection bin syst}. 

\subsection{Binomial systems}\label{subsection bin syst}
D.~Mall formulated his algorithm in a purely combinatorial language. So, we want to go into the combinatorial details of \cite{M2000}. We use the notations of section \ref{subsection notations}.

\begin{definition}\label{def - binomial system - def}
A triple $(A,C,\rho)$ consisting of two subsets $A,C\subset \N_d^n$ and of an $n$-tuple $\rho\in\Z^n$ is called a \emph{binomial system} (of degree $d$ in $n$ indeterminates) if the following conditions hold:
\begin{enumerate}
\item  $C+\rho\subset\N^n_d$,
\item $A\cap C=A\cap(C+\rho)=C\cap(C+\rho)=\emptyset$,
\item $A\cup C$ and $A\cup(C+\rho)$ are Borel sets.
\end{enumerate}
\end{definition}

\begin{remark}\label{remark on binomial systems}
If $(A,C,\rho)$ is a binomial system, we always assume that it is of degree $d$ in $n$ indeterminates unless otherwise stated.

If $(A,C,\rho)$ is a binomial system with $C\neq\emptyset$, then property (i) implies that $\rho_1+\dots+\rho_n=0$.

If $(A,C,\rho)$ is a binomial system, then $A$ is a Borel set.

If $(A,C,\rho)$ is a binomial system, then for any term order of $\N^n$ we have: \linebreak[4]If $c<c+\rho$ for some $c\in C$, then $c<c+\rho$ for all $c\in C$.

If $(A,C,\rho)$ is a binomial system, then $(A,C+\rho,-\rho)$ is also a binomial system. Hence if $C\neq\emptyset$, we always may assume that $\rho_{m(\rho)}>0$.
\end{remark}

\begin{notation}\label{not-binomial ideal-not}
If $C\subset\N^n$ and $\rho\in\Z^n$ are such that $C+\rho\subset\N^n$, set 
$$\Bin(C,\rho):=\{X^c-X^{c+\rho}\in S\mid c\in C\}.$$
If $(A,C,\rho)$ is a binomial system, define the ideal 
$$F(A,C,\rho):=(X^A\cup\Bin(C,\rho)).$$
\end{notation}

\begin{example}\label{example2}
Let $n:=4$, $d:=3$, $A:=\{(3, 0, 0, 0), (2, 1, 0, 0), (2, 0, 1, 0), (2, 0, 0, 1)$, $(1, 2, 0, 0), (1, 1, 1, 0), (1, 1, 0, 1), (1, 0, 2, 0), (0, 3, 0, 0), (0, 2, 1, 0)\}$, $C:=\{(0,2,0,1)\}$, $\rho:=(1,-2,1,0)$. Then $(A,C,\rho)$ is a binomial system. The saturation of the ideal 
$$F(A,C,\rho)=(x^3, x^2y, x^2z, x^2t, xy^2, xyz, xyt, xz^2, y^3, y^2z, y^2t - xzt)\subset K[x,y,z,t]$$ 
is the Mall ideal $\idealc=(y^2-xz, x^2, xy, xz^2)$ of Example \ref{example1}.
\end{example}

It is very easy to compute Gr\"obner bases of ideals generated by binomial systems and of their saturations:

\begin{proposition}\label{Groebner bases of (saturated) binomial ideals}
Let $(A,C,\rho)$ be a binomial system.

a) $X^A\cup\Bin(C,\rho)$ is a Gr\"obner basis of $F(A,C,\rho)$ with respect to the reverse lexicographic order. Especially, if $\rho_{m(\rho)}>0$, then $\init_\degrevlex F(A,C,\rho)=(X^{A\cup C})$.

b) $F(A,C,\rho)^\sat=(F(A,C,\rho):X_n^\infty)$.

c) If $\rho_{m(\rho)}>0$, then $X^{A^*}\cup\Bin(C^*,\rho)$ is a Gr\"obner basis of $F(A,C,\rho)^\sat$ with respect to the reverse lexicographic order.
\end{proposition}

\begin{proof}
The first two properties are shown in \cite[3.7(1)]{M2000} and \cite[3.8]{M2000}. Assertion c) follows from a), b) and the characteristic property of the reverse lexicographic order (cf. \cite[12.1]{Stu}).
\end{proof}

The following definition is crucial for section \ref{subsection generic initial ideal}.
\begin{definition}\label{def - admissible, good, Mall - def}
A binomial system $(A,C,\rho)$ is \emph{good} if $b_i=c_i$ for all $1\leq i<m(\rho)$ and for all $b,c\in C$.
\end{definition}

Mall did not state the following fact explicitly. It is proved in detail in section 2.1 of \cite{F}.
\begin{proposition}\label{Mall ideals are good}
Let $\idealc\subset S$ be a Mall ideal. Then there exists $d\in\N$ and a good binomial system $(A,C,\rho)$ of degree $d$ in $n$ indeterminates such that $\idealc=F(A,C,\rho)^\sat$.
\end{proposition}

\subsection{Sequentially Cohen-Macaulay\-ness}\label{subsection seq. CM}
In \cite[2.2]{HS} J.~Herzog and E.~Sbarra showed that in characteristic zero the Borel ideals are sequentially Cohen-Macaulay. Later J.~Herzog, D.~Popescu, and M.~Vladoiu \cite{HPV} generalized this result to monomial ideals of Borel type in any characteristic of $K$. An ideal $\idealb\subset S$ is of Borel type if $(\idealb:X_j^\infty)=(\idealb:_S\langle X_1,\dots,X_j\rangle_S^\infty)$ for all $1\leq j\leq n$. It is well known that Borel-fixed ideals are of Borel type (\cite[15.24]{E}). Since an ideal $F(A,C,\rho)$ generated by a good binomial system is fixed under the action of the unipotent group (cf.~Proposition \ref{good binomial ideals are fixed by unipotents}), it is not so surprising that $F(A,C,\rho)$ is sequentially Cohen-Macaulay. Indeed $F(A,C,\rho)$ is sequentially Cohen-Macaulay, even if $(A,C,\rho)$ is not good (cf.~Proposition \ref{binomial ideals and sCM}).

\begin{definition}
A homogeneous ideal $\ideala\subset S$ is \emph{sequentially Cohen-Macaulay} if there exists a finite filtration
$$\ideala=\ideala_0\subsetneq \ideala_1\subsetneq\dots\subsetneq \ideala_r=S$$ by homogeneous ideals such that
\begin{enumerate}
\item $\ideala_i/\ideala_{i-1}$ is Cohen-Macaulay for all $1\leq i\leq r$,
\item $\dim(\ideala_i/\ideala_{i-1})<\dim(\ideala_{i+1}/\ideala_i)$ for all $1\leq i<r$.
\end{enumerate}
\end{definition}

\sk
Given a binomial system $(A,C,\rho)$, we construct a filtration 
$$F(A,C,\rho)=F_0\subset F(A,C,\rho)^\sat=F_1\subset\dots\subset F_n=S$$ 
of ideals such that the quotients $F_{i+1}/F_i$ are zero or Cohen-Macaulay of dimension $i$ for all $0\leq i<n$ (Proposition \ref{Binomial quotients are CM}). There is a natural way to define the this filtration:
\begin{notation}
For $i\in\N$ let $S_{(i)}$ denote the polynomial ring $K[X_1,\dots,X_i]$ and let $\terms[i]$ be the set of monomials of $S_{(i)}$. In particular,  $S_{(n)}=S$ and $\terms[n]=\terms$. For $i<j$ one has a canonical inclusion $S_{(i)}\subset S_{(j)}$. If $i\in\N$ and $M\subset S_{(i)}$, we write $\langle M\rangle_{S_{(i)}}$ for the ideal in $S_{(i)}$ generated by $M$.

For $A\subset\N^n$ and $1\leq i\leq n$ set 
$$A_i:=\{(a_1,\dots,a_i)\in\N^i\mid(a_1,\dots,a_i,0,\dots,0)\in A\}.$$

Let $(A,C,\rho)$ be a binomial system. Set 
$$F_0:=F(A,C,\rho).$$
For $1\leq i\leq n-m(\rho)+1$ set
$$F_i:=(F(A,C,\rho)\cap S_{(n-i+1)})^\sat S.$$
For $n-m(\rho)+2\leq i\leq n$ set
$$F_i:=\bigl((\init_\degrevlex F(A,C,\rho))\cap S_{(n-i+1)}\bigr)^\sat S.$$
\end{notation}

It is easy to see that $F_0\subset\dots\subset F_{n-m(\rho)+1}$ and $F_{n-m(\rho)+2}\subset\dots\subset F_n$ (cf.~Lemma \ref{lemma about F_i(A,C,rho)}). The essential point is to show that $F_{n-m(\rho)+1}\subset F_{n-m(\rho)+2}$ and that $F_{n-m(\rho)+2}/F_{n-m(\rho)+1}$ is Cohen-Macaulay.

\sk
The following Lemma is crucial, because it allows us to compute the generators and the initial ideals of all ideals $F_i$.
\begin{lemma}\label{Groebner basis of F_i(A,C,rho)}
Let $(A,C,\rho)$ be a binomial system with $\rho_{m(\rho)}>0$, and let $m(\rho)\leq i\leq n$. Then 
$$X^{(A_i)^*}\cup\Bin((C_i)^*,(\rho_1,\dots,\rho_i))$$ 
is a Gr\"obner basis of $(F(A,C,\rho)\cap S_{(i)})^\sat$ with respect to the reverse lexicographic order of\/ $\terms[i]$.
\end{lemma}

\begin{proof}
Since $(A_i,C_i,(\rho_1,\dots,\rho_i))$ is a binomial system of degree $d$ in $i$ indeterminates, by Proposition \ref{Groebner bases of (saturated) binomial ideals}~c) the set $X^{(A_i)^*}\cup\Bin((C_i)^*,(\rho_1,\dots,\rho_i))$ is a Gr\"obner basis of $F(A_i,C_i,(\rho_1,\dots,\rho_i))^\sat$ with respect to the reverse lexicographic order of $\terms[i]$. To complete the proof it is enough to show that $$F(A,C,\rho)\cap S_{(i)}=F(A_i,C_i,(\rho_1,\dots,\rho_i)).$$ Since $m(\rho)\leq i$, the inclusion ``$\supset$'' is obvious. 

Let $f\in F(A,C,\rho)\cap S_{(i)}=\langle X^A\cup\Bin(C,\rho)\rangle_S\cap S_{(i)}$, and write 
$$f=\sum_{j=1}^r u_j X^{a_j}+\sum_{i=j}^s v_j (X^{c_j}-X^{c_j+\rho})$$ 
with $u_1$, \dots, $u_r$, $v_1$, \dots, $v_s\in\terms$ and $a_1$, \dots, $a_r\in A$, $c_1$, \dots, $c_s\in C$. It follows from $m(\rho)\leq i$, that $X^c\in S_{(i)}$ if and only if $X^{c+\rho}\in S_{(i)}$ for all $c\in C$. Since $f\in S_{(i)}$, we may assume that $X^{a_j}$, $X^{c_j}-X^{c_j+\rho}\in S_{(i)}$ for all $j$. This shows 
$$f\in\langle X^{A_i}\cup\Bin(C_i,(\rho_1,\dots,\rho_i))\rangle_{S_{(i)}}=F(A_i,C_i,(\rho_1,\dots,\rho_i)),$$ and the proof is finished.
\end{proof}

\begin{lemma}\label{lemma about F_i(A,C,rho)}
Let $(A,C,\rho)$ be a binomial system with $\rho_{m(\rho)}>0$, and put $m:=m(\rho)$. Then

a) $F_i=(F_i\cap S_{(n-i)})S\ \  \forall\,i\in\{0,\dots,n\}\setminus\{n-m+1\}$,

b) $F_{i+1}=(F_i\cap S_{(n-i)})^\sat S\ \  \forall\,i\in\{0,\dots,n-1\}\!\setminus\!\{n-m+1\}$,

c) $(F(A,C,\rho)\cap S_{(m)})^\sat\subset\bigl((\init_\degrevlex F(A,C,\rho))\cap S_{(m-1)}\bigr)^\sat S_{(m)},$

d) $F_i\subset F_{i+1}$ for all\/ $0\leq i<n$. 

e) $\{X_1,\dots,X_{m-1}\}\subset$
$$\sqrt{\Bigl((F(A,C,\rho)\cap S_{(m)})^\sat\smash{\underset{S_{(m)}}{:}}\bigl((\init_\degrevlex F(A,C,\rho))\cap S_{(m-1)}\bigr)^\sat S_{(m)}\Bigr)}.$$

f) Let $u\in\terms[m]$. Then 
$u\notin\init_\degrevlex(F(A,C,\rho)\cap S_{(m)})^\sat$ implies 
$$X_m u\notin\init_\degrevlex(F(A,C,\rho)\cap S_{(m)})^\sat.$$
\end{lemma}

\smallskip
\begin{proof}
a) Let $i\in\{0,\dots,n\}\setminus\{n-m+1\}$. If $i>n-m+1$, then by Remark \ref{Groebner bases of Borel ideals} the ideal $F_i=\bigl((\init_\degrevlex F(A,C,\rho))\cap S_{(n-i+1)}\bigr)^\sat S$ is generated in $S_{(n-i)}$ since $(\init_\degrevlex F(A,C,\rho))\cap S_{(n-i+1)}\subset S_{(n-i+1)}$ is a Borel ideal by Proposition \ref{Groebner bases of (saturated) binomial ideals}~a). 

If $n-i+1>m$, then by Lemma \ref{Groebner basis of F_i(A,C,rho)} we have $F_i=(F(A,C,\rho)\cap S_{(n-i+1)})^\sat S=\bigl\langle X^{(A_{n-i+1})^*}\cup\Bin((C_{n-i+1})^*,(\rho_1,\dots,\rho_{n-i+1}))\bigr\rangle_{S_{(n-i+1)}} S.$ Since $\rho_{n-i+1}=0$, the latter ideal is also generated in $S_{(n-i)}$.

\sk
b) Let $i\in\{0,\dots,n-1\}\setminus\{n-m+1\}$. Set 
$$\ideala:=
\begin{cases}
F(A,C,\rho),&\text{if $i<n-m+1$;}\\
\init_\degrevlex F(A,C,\rho),&\text{otherwise.}
\end{cases}$$
Since $i\neq n-m+1$ we have $F_i=(\ideala\cap S_{(n-i+1)})^\sat S$ and $F_{i+1}=(\ideala\cap S_{(n-i)})^\sat S$. It follows that
\begin{align*}
F_{i+1}&=(\ideala\cap S_{(n-i)})^\sat S=((\ideala\cap S_{(n-i+1)})\cap S_{(n-i)})^\sat S\\
&=((\ideala\cap S_{(n-i+1)})^\sat\cap S_{(n-i)})^\sat S=(F_i\cap S_{(n-i)})^\sat S.
\end{align*}

\sk
c) Set $\idealb:=\bigl((\init_\degrevlex F(A,C,\rho))\cap S_{(m-1)}\bigr)^\sat S_{(m)}.$ We first prove the following claim 
$$X^{(A_m)^*}\cup X^{(C_m)^*}\cup X^{(C_m)^*+(\rho_1,\dots,\rho_m)}\subset\idealb.$$
Let $b\in(A_m)^*\cup(C_m)^*\cup ((C_m)^*+(\rho_1,\dots,\rho_m))$. 

\emph{Case 1:} $b\in(A_m)^*\cup(C_m)^*$. Then there exists $a\in A\cup C$ such that 
$b=(a_1,\dots,a_{m-1},0)$ and $m(a)\leq m$. Let $a':=a+a_m(e_{m-1}-e_m)$. Since $A\cup C$ is a Borel set, we have $a'\in A\cup C$. Hence $(a_1,\dots,a_{m-2},a_{m-1}+\nolinebreak a_m)\in (A\cup C)_{m-1}$ and $(a_1,\dots,a_{m-2},0)\in ((A\cup C)_{m-1})^*$. By Remark~\ref{Groebner bases of Borel ideals} and Proposition \ref{Groebner bases of (saturated) binomial ideals}~a) it follows that 
$X^{(a_1,\dots,a_{m-2},0)}\in\langle X^{(A\cup C)_{m-1}}\rangle_{S_{(m-1)}}^\sat\subset\bigl((\init_\degrevlex F(A,C,\rho))\cap S_{(m-1)}\bigr)^\sat,$ and therefore, $X^b=X_{m-1}^{a_{m-1}}X^{(a_1,\dots,a_{m-2},0,0)}\in\idealb.$

\emph{Case 2:} $b\in(C_m)^*+(\rho_1,\dots,\rho_m)$. Then there exists $c\in C$ such that $b=(c_1,\dots,c_{m-1},0)+(\rho_1,\dots,\rho_m)$. By similar arguments to those used in case 1 we obtain $a':=c+\rho+(c+\rho)_m(e_{m-1}-e_m)\in A\cup(C+\rho)$. Now, $m(a')=m-1$ implies $a'\in A$, and therefore 
$X^{(c_1,\dots,c_{m-2},0)+(\rho_1,\dots,\rho_{m-2},0)}\in\bigl((\init_\degrevlex F(A,C,\rho))\cap S_{(m-1)}\bigr)^\sat.$ 
It follows again that $X^b=X_{m-1}^{(c+\rho)_{m-1}}X_m^{\rho_m}X^{(c_1,\dots,c_{m-2},0,0)+(\rho_1,\dots,\rho_{m-2},0,0)} \in\idealb,$ and our claim is proved.

Now we get our statement by means of Lemma \ref{Groebner basis of F_i(A,C,rho)}:
\begin{align*}
(F(A,C,\rho)\cap S_{(m)})^\sat&=\langle X^{(A_m)^*}\cup\Bin((C_m)^*,(\rho_1,\dots,\rho_m))\rangle_{S_{(m)}}\\*
&\subset\langle X^{(A_m)^*}\cup X^{(C_m)^*}\cup X^{(C_m)^*+(\rho_1,\dots,\rho_m)}\rangle_{S_{(m)}}\subset\idealb.
\end{align*}

\sk
d) This statement follows immediately from a), b), and c).

\sk
e) By Proposition \ref{Groebner bases of (saturated) binomial ideals}~a) and Remark \ref{Groebner bases of Borel ideals} it holds 
\begin{align*}
\bigl((\init_\degrevlex F(A,C,\rho))\cap S_{(m-1)}\bigr)^\sat&=(\langle X^{A\cup C}\rangle_S\cap S_{(m-1)})^\sat\\*
&=\langle X^{((A\cup C)_{m-1})^*}\rangle_{S_{(m-1)}}.
\end{align*}
Without loss of generality we may assume that $A_{m-1}\neq\emptyset$ and set 
$$r:=\max{\{a_{m-1}\in\N\mid a\in (A\cup C)_{m-1}\}}+\rho_m.$$ Let $1\leq i<m$. It is enough to show that 
$$X_i^r\langle X^{((A\cup C)_{m-1})^*}\rangle_{S_{(m-1)}}S_{(m)}\subset(F(A,C,\rho)\cap S_{(m)})^\sat.$$
Let $a\in(A\cup C)_{m-1}$. Since $(A\cup C)_{m-1}$ is a Borel set, we have 
$$c:=a_{m-1}e_i+a^*=a+a_{m-1}(e_i-e_{m-1})\in(A\cup C)_{m-1}.$$ 

\emph{Case 1:} $c\in A_{m-1}$. Then 
$$X_i^{a_{m-1}}X^{a^*}\in F(A,C,\rho)\cap S_{(m-1)}\subset(F(A,C,\rho)\cap S_{(m)})^\sat.$$

\emph{Case 2:} $c\in C_{m-1}$. Then $c':=(c_1,\dots,c_{m-1},0,\dots,0)\in C$. Since $(A\cup(C+\rho))$ is a Borel set, we have $b:=c'+\rho+\rho_m(e_i-e_{m})\in A\cup(C+\rho).$ It is clear that $m(b)<m$, whence $b\in A$. It follows that $X_i^{a_{m-1}+\rho_m}X^{a^*}=x_i^{\rho_m}(X^{c'}-X^{c'+\rho})+X_m^{\rho_m}X^c\in F(A,C,\rho)\cap S_{(m)}\subset(F(A,C,\rho)\cap S_{(m)})^\sat$.

In both cases our claim follows.

\sk
f) From Lemma \ref{Groebner basis of F_i(A,C,rho)} it follows that 
$$\init_\degrevlex(F(A,C,\rho)\cap S_{(m)})^\sat=\langle X^{(A_m\cup C_m)^*}\rangle_{S_{(m)}}.$$ 
Thus, if $X_m u\in\init_\degrevlex(F(A,C,\rho)\cap S_{(m)})^\sat$, there exists $b\in(A_m\cup C_m)^*$ such that $X^b$ divides $X_m u$. Since $X^b$ is not divisible by $X_m$, it has to divide $u$. Hence  $u\in\init_\degrevlex(F(A,C,\rho)\cap S_{(m)})^\sat$.
\end{proof}

\begin{proposition}\label{Binomial quotients are CM}
Let $(A,C,\rho)$ be a binomial system, and let\/ $0\leq i<n$. Then $F_{i+1}/F_i$ is zero or Cohen-Macaulay of dimension $i$.
\end{proposition}

\begin{proof}
Let $m:=m(\rho)$. By Remark \ref{remark on binomial systems} we may assume that $\rho_m>0$. Assume first that $i\neq n-m+1$. Set $\ideala:=F_i\cap S_{(n-i)}.$ Then by Lemma \ref{lemma about F_i(A,C,rho)}~a) and~b) we have $F_i=\ideala\,S$ and $F_{i+1}=\ideala^\sat S$. It follows that
$$F_{i+1}/F_i=\ideala^\sat S/\ideala\,S=(\ideala^\sat/\ideala)\otimes_{S_{(n-i)}}S=H^1_{(S_{(n-i)})_+}(\ideala)\otimes_{S_{(n-i)}}S.$$
The (finitely generated) $S_{(n-i)}$-module $H^1_{(S_{(n-i)})_+}\!(\ideala)$ is Artinian, therefore it is zero or Cohen-Macaulay and zero-dimensional. Thus, the $S$-module $F_{i+1}/F_i$ is zero or Cohen-Macaulay of dimension $i$ (cf.\ \cite[2.1.9]{BH}).

\sk
We now prove our statement for $i=n-m+1$. Set
$$\ideala:=\bigl((\init_\degrevlex F(A,C,\rho))\cap S_{(m-1)}\bigr)^\sat S_{(m)},$$
$$\idealb:=(F(A,C,\rho)\cap S_{(m)})^\sat.$$
By Lemma \ref{lemma about F_i(A,C,rho)}~c) it holds $\idealb\subset\ideala$.
Since 
\begin{align*}
F_{n-m+2}/F_{n-m+1}&=\bigl((\init_\degrevlex F(A,C,\rho))\cap S_{(m-1)}\bigr)^\sat S/(F(A,C,\rho)\cap S_{(m)})^\sat S\\*
&=(\ideala/\idealb)\otimes_{S_{(m)}}S,
\end{align*}
it is enough to show that the $S_{(m)}$-module $\ideala/\idealb$ is zero or one-dimensional and Cohen-Macaulay.

By Lemma \ref{lemma about F_i(A,C,rho)}~e) it holds
\begin{align*}
\dim \ideala/\idealb&=\dim S_{(m)}/(\idealb\underset{S_{(m)}}{:}\ideala)=\dim S_{(m)}/\sqrt{(\idealb\smash{\underset{S_{(m)}}{:}}\ideala)}\\
&\leq \dim S_{(m)}/\langle X_1,\dots,X_{m-1}\rangle_{S_{(m)}}=1.
\end{align*}
Thus, it is enough to show that $X_m$ is a non-zerodivisor of $\ideala/\idealb$. 

Lemma \ref{lemma about F_i(A,C,rho)}~f) states that $X_m$ is a non-zerodivisor of $S_{(m)}/\init_\degrevlex\idealb$. Since any non-zerodivisor of $S_{(m)}/\init_\degrevlex\idealb$ is a non-zerodivisor of $S_{(m)}/\idealb$ and hence of $\ideala/\idealb$, our proof is finished.
\end{proof}

Cancelling the redundant ideals in the filtration $F_0\subset F_1\subset\dots\subset F_n=S$ we get the following statement:

\begin{proposition}\label{binomial ideals and sCM}
Let $(A,C,\rho)$ be a binomial system. Then the ideals $F(A,C,\rho)$ and $F(A,C,\rho)^\sat$ are sequentially Cohen-Macaulay.
\end{proposition}

Together with Proposition \ref{Mall ideals are good} we have proved:
\begin{theorem}\label{Mall ideals are sCM}
Mall ideals are sequentially Cohen-Macaulay.
\end{theorem}

\subsection{The generic initial ideal of a binomial ideal}\label{subsection generic initial ideal}
Let $(A,C,\rho)$ be a good binomial system. Then the generic initial ideal of $F(A,C,\rho)$ and the initial ideal of $F(A,C,\rho)$ with respect to any admissible term order coincide (Proposition \ref{Gin = in for good binomial ideals}). This conclusion is not at all trivial: It does not hold if $(A,C,\rho)$ is not good (s.~Example \ref{counterexample}). The crucial point is that $F(A,C,\rho)$ is fixed under the action of the unipotent group, if $(A,C,\rho)$ is good (Proposition \ref{good binomial ideals are fixed by unipotents}). In the following we will prove this statement.

In case the binomial system  $(A,C,\rho)$ is good, we want to compute $g(f)$ if $f$ is a generator of the ideal $F(A,C,\rho)$ and $g\in\unipotent$ is unipotent. To do this, we introduce generic coordinates for $g$: 
\begin{notation}\label{not generic coordinates}
Set $T:=K[Y_{ij}\mid 1\leq i\leq j\leq n]$.

Define the automorphism of $T$-algebras 
$$\phi:T[X_1,\dots,X_n]\ra T[X_1,\dots,X_n]$$ 
by $\phi(X_j):=\sum_{i=1}^j Y_{ij}X_i$.

Let $g=[g_{ij}\mid 1\leq i\leq n,\:1\leq j\leq n]\in K^{(n,n)}$ be a matrix. Then we denote by $\bar{g}:T[X_1,\dots,X_n]\ra S$ the homomorphism of $S$-algebras defined by $\bar{g}(Y_{ij}):=g_{ij}$ for $1\leq i\leq j\leq n$.

For $a\in\N^n$ set $g^a:=\prod_{i=1}^n g_{ii}^{a_i}$ and $Y^a:=\prod_{i=1}^n Y_{ii}^{a_i}$.

For $M\in\N^{(n,n)}$ set $Y^M:=\prod_{1\leq i\leq j\leq n}Y_{ij}^{M_{ij}}$.

For $\rho\in\Z^n$ and $M\in\Z^{(n,n)}$ let $M+\rho\in\Z^{(n,n)}$ denote the matrix which is defined by
$$(M+\rho)_{ij}:=
\begin{cases}
M_{ij},&\text{if $i\neq j$;}\\
M_{ii}+\rho_i,&\text{otherwise.}
\end{cases}$$

For $a,b\in\N^n$ set 
$$U(a,b):=\{M\in U(n)\mid \sum_{i=1}^n M_{ji}=a_j,\sum_{i=1}^n M_{ij}=b_j\ \forall\ 1\leq j\leq n\},$$
where $U(n)$ is defined as in \ref{combinatorial notations}.
\end{notation}

\begin{remark}\label{remark on U(a,b)}
Let $a,b\in\N^n$. Then by Lemma \ref{properties of Borel order}: $a\borelgeq b\iff U(a,b)\neq\emptyset$.
\end{remark}

\begin{lemma}\label{lemma on Borel matrices}
Let $b,c\in\N^n_d$ and $\rho\in\Z^n$ be such that $b+\rho$, $c+\rho\in\N^n_d$ and such that $b_i=c_i$ for all $1\leq i<m(\rho)$.

a) Let $b\borelgeq c$ and $M\in U(b,c)$. Then $M_{jj}=c_j$ for $1\leq j\leq m(\rho)$.

b) Let $b\borelgeq c$ and $M\in U(b,c)$. Then $M_{jj}+\rho_j\geq 0$ for $1\leq j\leq n$.

c) Let $M\in\Z^{(n,n)}$. Then $M\in U(b,c)\iff M+\rho\in U(b+\rho,c+\rho)$.
\end{lemma}

\begin{proof}
a) We use an inductive argument. It is clear that $M_{11}=\sum_{i=1}^n M_{i1}=c_1$. Let $1<j\leq m(\rho)$ and assume that $c_k=M_{kk}$ for all $1\leq k<j$. We have $M_{kk}=c_k=b_k=\sum_{i=1}^n M_{ki}$ for $1\leq k<j$ and therefore $M_{kj}=0$ for $1\leq k<j$. Since $M_{ij}=0$ for all $i>j$, it holds $c_j=\sum_{i=1}^n M_{ij}=M_{jj}$.

\sk
b) This statement follows from statement a) and from the condition $c_j+\rho_j\geq 0$ for $1\leq j\leq n$.

\sk
c) Let $M\in U(b,c)$. Then part b) states that $M+\rho\in U(n)$. Since
$$(b+\rho)_j=b_j+\rho_j=\sum_{i=1}^n M_{ji}+\rho_j=\sum_{i=1}^n(M+\rho)_{ji}$$
and
$$(c+\rho)_j=c_j+\rho_j=\sum_{i=1}^n M_{ij}+\rho_j=\sum_{i=1}^n(M+\rho)_{ij}$$
for all $1\leq j\leq n$, it holds $M+\rho\in U(b+\rho,c+\rho)$.

In order to prove the converse implication, one just has to replace $b$, $c$, $\rho$ with $b+\rho$, $c+\rho$, $-\rho$ respectively.
\end{proof}

\begin{notation}\label{not coefficients}
For $m$ and $m_1,\dots,m_n\in\N$ with $m=\sum_{i=1}^n m_i$ set 
$$\binom{m}{m_1,\dots,m_n}:=\frac{m!}{m_1!\cdots m_n!}.$$

For $M\in\N^{(n,n)}$ set
$$\mu_M:=\prod_{j=1}^n\binom{\sum_{i=1}^n M_{ij}}{M_{1j},\dots,M_{nj}}.$$

For $a,b\in\N^n_d$ let $\alpha_a^b\in T$ be the coefficient of $X^a$ in the polynomial $\phi(X^b)$, so that $\phi(X^b)=\sum_{a\in\N^n_d}\alpha^b_a X^a$.

For $b,c\in\N^n_d$ and $\rho\in\Z^n$ such that $b+\rho$, $c+\rho\in\N^n_d$ and such $b_i=c_i$ for all $1\leq i<m(\rho)$ set 
$$p^\rho_{b,c}:=\sum_{M\in U(b,c)}\mu_M Y^{M-\rho^-},$$
where $\rho^-$ is defined according to \ref{combinatorial notations}.
\end{notation}

\begin{remark}\label{remark on p^rho}
Let $b,c\in\N^n_d$ and $\rho\in\Z^n$ be such that $b+\rho$, $c+\rho\in\N^n_d$ and such $b_i=c_i$ for all $1\leq i<m(\rho)$. Then for all $M\in U(b,c)$ we have by Lemma \ref{lemma on Borel matrices} b) that $M_{jj}+\rho_j\geq 0$ for $1\leq j\leq n$, whence $Y^{M-\rho^-}\in T[X_1,\dots,X_n]$. Therefore, $p^\rho_{b,c}$ is a polynomial in $T[X_1,\dots,X_n]$.

If $b\not\borelgeq c$, then $p^\rho_{b,c}=0$ by Remark \ref{remark on U(a,b)}.
\end{remark}

\begin{lemma}\label{computation of alpha^b_a}
Let $a,b\in\N^n_d$. Then $\alpha^b_a=\sum_{M\in U(a,b)}\mu_M Y^M$. In particular $a\borelgeq b$ if and only if $\alpha^b_a\neq 0$.
\end{lemma}

\begin{proof}
From the multinomial formula it follows
\begin{align*}
\sum_{a\in\N^n_d}\alpha^b_a X^a&=\phi(X^b)=\prod_{j=1}^n \bigl(\sum_{i=1}^j Y_{ij}X_i\bigr)^{b_j}\\
&=\prod_{j=1}^n\sum_{\substack{k_1,\,\dots,\,k_j\in\N\\k_1+\dots+k_j=b_j}}\binom{b_j}{k_1,\dots,k_j}\prod_{i=1}^j (Y_{ij}X_i)^{k_i}\\
&=\sum_{\substack{k_{11},\,k_{12},\,k_{22},\,\dots,\,k_{1n},\,\dots,\,k_{nn}\in\N\\\sum_{i=1}^j k_{ij}=b_j\ \forall\ 1\leq j\leq n}}\ \prod_{j=1}^n\binom{b_j}{k_{1j},\dots,k_{jj}}\prod_{i=1}^j (Y_{ij}X_i)^{k_{ij}}\\
&=\sum_{\substack{M\in U(n)\\\sum_{i=1}^n M_{ij}=b_j\ \forall\ 1\leq j\leq n}}\prod_{j=1}^n\binom{b_j}{M_{1j},\dots,M_{nj}}\prod_{i=1}^j (Y_{ij}X_i)^{M_{ij}}\\
&=\sum_{\substack{M\in U(n)\\\sum_{i=1}^n M_{ij}=b_j\ \forall\ 1\leq j\leq n}}\mu_M\prod_{1\leq i\leq j\leq n} Y_{ij}^{M_{ij}}X_i^{M_{ij}}\\
&=\sum_{\substack{M\in U(n)\\\sum_{i=1}^n M_{ij}=b_j\ \forall\ 1\leq j\leq n}}\mu_M Y^M\prod_{i=1}^n X_i^{\sum_{j=1}^n M_{ij}}\\
&=\sum_{\substack{M\in U(n)\\\sum_{i=1}^n M_{ij}=b_j\ \forall\ 1\leq j\leq n}}\mu_M Y^M X^{(\sum_{j=1}^n M_{1j},\dots,\sum_{j=1}^n M_{nj})}.
\end{align*}
Hence we get 
$$\alpha^b_a=\sum_{\substack{M\in U(n)\\\sum_{i=1}^n M_{ij}=b_j\ \forall\ 1\leq j\leq n\\\sum_{j=1}^n M_{ij}=a_i\ \forall\ 1\leq i\leq n}} \mu_M Y^M=\sum_{M\in U(a,b)}\mu_M Y^M.$$

In particular, $a\borelgeq b\iff U(a,b)\neq\emptyset\iff\alpha^b_a\neq 0$, by Remark \ref{remark on U(a,b)}.
\end{proof}

\begin{lemma}\label{computation of alpha^c_b}
Let $b,c\in\N^n_d$ and $\rho\in\Z^n$ be such that $b+\rho$, $c+\rho\in\N^n_d$ and $b_i=c_i$ for all $1\leq i<m(\rho)$. Then  $\alpha^c_b=p^\rho_{b,c}\,Y^{\rho^-}$\! and $\alpha^{c+\rho}_{b+\rho}=p^\rho_{b,c}\,Y^{\rho^+}$.
\end{lemma}

\begin{proof}
If $b\not\borelgeq c$, then $\alpha^c_b=\alpha^{c+\rho}_{b+\rho}=0$ by Lemma \ref{computation of alpha^b_a} and $p^\rho_{b,c}=0$ by Remark \ref{remark on p^rho}. Hence, we can assume that $b\borelgeq c$.
The first equation follows immediately from Lemma \ref{computation of alpha^b_a}. To prove the second one, we first claim that $\mu_M=\mu_{M+\rho}$ for all $M\in U(b,c)$. Let $M\in U(b,c)$. Let $1\leq j\leq n$. If $j\leq m(\rho)$ we have by part a) of Lemma \ref{lemma on Borel matrices} that
$$\binom{c_j}{M_{1j},\dots,M_{nj}}=\frac{c_j!}{M_{jj}!}=1=\frac{(c_j+\rho_j)!}{(M_{jj}+\rho_j)!}=\binom{c_j+\rho_j}{(M+\rho)_{1j},\dots,(M+\rho)_{nj}}.$$
If $j>m(\rho)$ we have $\binom{c_j}{M_{1j},\dots,M_{nj}}=\binom{c_j+\rho_j}{(M+\rho)_{1j},\dots,(M+\rho)_{nj}}$. Our claim now follows from the definition of $\mu_M$ and $\mu_{M+\rho}$.

By Lemma \ref{computation of alpha^b_a} and \ref{lemma on Borel matrices}~c) we get
$$\alpha ^{c+\rho}_{b+\rho}=\sum_{M\in U(b+\rho,c+\rho)}\mu_{M+\rho}Y^{M+\rho}=\sum_{M\in U(b,c)}\mu_M Y^{M-\rho^-+\rho^+}=p^\rho_{b,c}\,Y^{\rho^+}.$$
\end{proof}

\begin{proposition}\label{good binomial ideals are fixed by unipotents}
Let $(A,C,\rho)$ be a good binomial system and $g\in\unipotent$ an unipotent matrix. Then $g(F(A,C,\rho))=F(A,C,\rho).$
\end{proposition}

\begin{proof}
If $C=\emptyset$, the statement follows from the fact that the ideal $(X^A)$ is Borel-fixed (s.~ Proposition~\ref{existence of gin}). 

Now let $C\neq\emptyset$ and $c:=\min_\degrevlex C$. Then $(A,C\setminus\{c\},\rho)$ is again a good binomial system (cf.\ Lemma \ref{Borel greater implies greater for any term order}). Since $g(F(A,C,\rho))=(g(X^A\cup\Bin(C,\rho))$, it is enough to show the following equality of $K$-vectorspaces:
$$\langle g(X^A\cup\Bin(C,\rho))\rangle_K=\langle X^A\cup\Bin(C,\rho)\rangle_K.$$ 
Hence it is enough to show that 
$$g(X^A\cup\Bin(C,\rho))\subset\langle X^A\cup\Bin(C,\rho)\rangle_K.$$
By induction on $\#C$ we may assume that $g(F(A,{C\setminus\{c\}},\rho))=F(A,{C\setminus\{c\}},\rho)$, whence 
$$g(X^A\cup\Bin(C\setminus\{c\},\rho))\subset\langle X^A\cup\Bin(C,\rho)\rangle_K.$$ 
Therefore, it is enough to show that 
$$g(X^c-X^{c+\rho})\in\langle X^A\cup\Bin(C,\rho)\rangle_K.$$

In the following we make use of the Notations \ref{not generic coordinates} and \ref{not coefficients}. Observe that $\bar g(Y^{\rho^-})=g^{\rho^-}=1$ and $\bar g(Y^{\rho^+})=g^{\rho^+}=1$ as $g$ is unipotent.
Since $(A,C,\rho)$ is good, we may apply Lemma \ref{computation of alpha^c_b} to compute
\begin{align*}
g(X^c-X^{c+\rho})&=\ \bar{g}(\phi(X^c-X^{c+\rho}))\\
&=\ \bar{g}\bigl(\sum_{a\in\N^n_d}\alpha^c_a X^a - \sum_{a\in\N^n_d}\alpha^{c+\rho}_a X^a\bigr)\\
&=\ \bar{g}\bigl(\sum_{a\borelgeq c}\alpha^c_a X^a - \hspace{-8pt}\sum_{a\borelgeq c+\rho}\hspace{-8pt} \alpha^{c+\rho}_a X^a\bigr)\\
&=\ \bar{g}\bigl(\sum_{\substack{a\in A\\a\borelg c}}\alpha^c_a X^a -\hspace{-8pt}\sum_{\substack{a\in A\\a\borelg c+\rho}}\hspace{-8pt} \alpha^{c+\rho}_a X^a\bigr)+\bar{g}\bigl(\!\sum_{\substack{a\in C\\a\borelgeq c}}\!\alpha^c_a X^a - \hspace{-8pt}\sum_{\substack{a\in C+\rho\\a\borelgeq c+\rho}}\hspace{-8pt} \alpha^{c+\rho}_a X^a\bigr)\\
&\hspace{-0.38pc}\overset{(\ref{linearity of Borel order})}=\,\bar{g}\bigl(\sum_{a\in A}\alpha^c_a X^a - \sum_{a\in A}\alpha^{c+\rho}_a X^a\bigr)+\bar{g}\bigl(\sum_{b\in C}(\alpha^c_b X^b - \alpha^{c+\rho}_{b+\rho} X^{b+\rho})\bigr)\\
&\hspace{-0.55pc}\overset{(\ref{computation of alpha^c_b})}{=}\bar{g}\bigl(\sum_{a\in A}(\alpha^c_a-\alpha^{c+\rho}_a)X^a\bigr)+\bar{g}\bigl(\sum_{b\in C} (p^\rho_{b,c}\,Y^{\rho^-}X^b -p^\rho_{b,c}\,Y^{\rho^+} X^{b+\rho})\bigr)\\
&=\ \sum_{a\in A}\bar{g}(\alpha^c_a-\alpha^{c+\rho}_a)X^a+\sum_{b\in C} \bar{g}(p^\rho_{b,c})\bigl(\bar{g}(Y^{\rho^-})X^b -\bar{g}(Y^{\rho^+}) X^{b+\rho}\bigr)\\
&=\ \sum_{a\in A}\bar{g}(\alpha^c_a-\alpha^{c+\rho}_a)X^a+ \sum_{b\in C} \bar{g}(p^\rho_{b,c})(X^b - X^{b+\rho}).
\end{align*}
Hence our Proposition is proved.
\end{proof}

Now, it follows immediately from Corollary \ref{Gin=in for unipotent fixed ideals}:

\begin{proposition}\label{Gin = in for good binomial ideals}
Let $(A,C,\rho)$ be a good binomial system and $\tau$ an admissible term order. Then $\Gin_\tau F(A,C,\rho)=\init_\tau F(A,C,\rho).$\hfill$\square$
\end{proposition}

As a further consequence we have by Proposition \ref{Mall ideals are good} and Proposition \ref{Gin and sat commute}:
\begin{theorem}\label{Gin = in for Mall ideals}
If $\idealc\subset S$ is a Mall ideal, then $\Gin_\degrevlex\idealc=\init_\degrevlex\idealc$.
\end{theorem}

\pagebreak[2]
We conclude this section with a example which shows that Proposition \ref{Gin = in for good binomial ideals} fails if the binomial system $(A,C,\rho)$ is not good:
\begin{example}\label{counterexample}
Consider the ring $R:=K[X_1,\dots,X_5]=K[x,y,z,t,u]$. Let $\rho:=(1,-2,2,-2,1)$,  $b:=(0,2,0,3,0)$, $c:=(0,2,0,2,1)$ and $C:=\{b,c\}$. Let $B\subset\N^5_5$ be the smallest Borel set containing $D:=C\cup(C+\rho)$ and set $A:=B\setminus D$. Then $(A,C,\rho)$ is a binomial system. Let $\ideala:=F(A,C,\rho)$. We then compute:\pagebreak[1]
\begin{align*}
\Gin_\degrevlex\ideala=(&
x^5, x^4y, x^3y^2, x^2y^3, xy^4, y^5, x^4z, x^3yz, x^2y^2z, xy^3z, y^4z, x^3z^2,\\* &x^2yz^2, xy^2z^2, y^3z^2, x^2z^3, xyz^3, y^2z^3, xz^4, x^4t, x^3yt, x^2y^2t,\\* &xy^3t, y^4t, x^3zt, x^2yzt, xy^2zt, y^3zt, x^2z^2t, xyz^2t, y^2z^2t, xz^3t,\\* &x^3t^2, x^2yt^2, xy^2t^2, y^3t^2, x^2zt^2, xyzt^2, y^2zt^2, xz^2t^2, x^2t^3, xyt^3, y^2t^3,\\* &x^4u, x^3yu, x^2y^2u, xy^3u, y^4u, x^3zu, x^2yzu, xy^2zu, y^3zu, x^2z^2u,\\* & xyz^2u, y^2z^2u, xz^3u, x^3tu, x^2ytu, xy^2tu, y^3tu, x^2ztu, xyztu,\\* &y^2ztu, \underline{xz^2tu}, x^2t^2u, xyt^2u, x^3u^2, x^2yu^2, xy^2u^2, x^2zu^2, xyzu^2),\\
\init_\degrevlex\ideala=(&
x^5, x^4y, x^3y^2, x^2y^3, xy^4, y^5, x^4z, x^3yz, x^2y^2z, xy^3z, y^4z, x^3z^2,\\* &x^2yz^2, xy^2z^2, y^3z^2, x^2z^3, xyz^3, y^2z^3, xz^4, x^4t, x^3yt, x^2y^2t,\\* &xy^3t, y^4t, x^3zt, x^2yzt, xy^2zt, y^3zt, x^2z^2t, xyz^2t, y^2z^2t, xz^3t,\\* &x^3t^2, x^2yt^2, xy^2t^2, y^3t^2, x^2zt^2, xyzt^2, y^2zt^2, xz^2t^2, x^2t^3, xyt^3, y^2t^3,\\* &x^4u, x^3yu, x^2y^2u, xy^3u, y^4u, x^3zu, x^2yzu, xy^2zu, y^3zu, x^2z^2u,\\ &xyz^2u, y^2z^2u, xz^3u, x^3tu, x^2ytu, xy^2tu, y^3tu, x^2ztu, xyztu,\\* &y^2ztu, x^2t^2u, xyt^2u, \underline{y^2t^2u}, x^3u^2, x^2yu^2, xy^2u^2, x^2zu^2, xyzu^2).
\end{align*}
These ideals are not equal, as is indicated by the underlined generators. The reason is the following: Let $M\in\N^{(5,5)}$ be the unique element of $U(b,c)$. Then $M_{55}=0\neq c_5$, but $m(\rho)=5$ (counterexample to part (a) of Lemma~\ref{lemma on Borel matrices}). We cannot conclude that $\mu_M$ equals $\mu_{M+\rho}$; indeed $\mu_M=1$ and $\mu_{M+\rho}=2$. It follows that $\alpha_{b+\rho}^{c+\rho}=2p^\rho_{b,c}Y^{\rho^+}$ (counterexample to Lemma~\ref{computation of alpha^c_b}). Let $g\in\Gl(5,K)$ be unipotent. We then compute 
\begin{multline*}
g(X^c-X^{c+\rho})=X^c - X^{c+\rho} + \bar{g}(p^\rho_{b,c})(X^b -2X^{b+\rho})+ \sum_{a\in A}\bar{g}(\alpha^c_a-\alpha^{c+\rho}_a)X^a\\*
\notin\langle X^A\cup\Bin(C,\rho)\rangle_K
\end{multline*}
(counterexample to Proposition \ref{good binomial ideals are fixed by unipotents}). A further computation yields 
$$g(X^b-X^{b+\rho})=X^b - X^{b+\rho} + \sum_{a\in A}\bar{g}(\alpha^b_a-\alpha^{b+\rho}_a)X^a.$$ 
Since $A$ is a Borel set, it is clear that $X^A\subset g(F(A,C,\rho))$, whence $X^b-X^{b+\rho}$, $X^b -2X^{b+\rho} + \bar{g}(p^\rho_{b,c})^{-1}(X^c - X^{c+\rho})\in g(F(A,C,\rho))$. Since $b>_\degrevlex c$ we get $X^{b+\rho}\in\Gin_\degrevlex F(A,C,\rho)\setminus\init_\degrevlex F(A,C,\rho)$ and $X^c\in\init_\degrevlex F(A,C,\rho)\setminus\Gin_\degrevlex F(A,C,\rho)$.
\end{example}

\sk
\section{Application to Hilbert function strata}\label{Hilbert scheme}
Let $p\in\Q[t]$ be a polynomial. The Hilbert scheme $\Hilb^p=\Hilb^p_{\proj}$ is defined to be the representing scheme of the Hilbert functor $\Hilbf^p_{\proj}:\sch_K\ra\sets$ which assigns to each locally Noetherian scheme $T$ over $K$ the set of all closed subschemes $W\subset\proj\times T$, flat over $T$ such that for every point $x\in T$ the fibre $W_x$ of $W$ over $x$ has Hilbert polynomial $p$ (cf. \cite{Gr}, \cite{St}). The Hilbert scheme is characterized by the following 

\begin{property}
There exists a universal closed subscheme $W_{\Hilb}\subset\proj\times\Hilb^p$, flat over $\Hilb^p$ with Hilbert polynomial $p$ in all fibres such that for every locally Notherian $K$-scheme $T$ and for every closed subschemes $W\subset\proj\times T$, flat over $T$ with Hilbert polynomial $p$ in all fibres there exists a unique morphism of $K$-schemes $g:T\ra\Hilb^p$ such that $W=W_{\Hilb}\times_{\Hilb^p}T$ is the pullback of the universal subscheme by $g$.
\end{property}

For each point $x\in\Hilb^p_K$ let $\kappa(x)$ denote the residue field of $x$ on $\Hilb^p$ and set $W_x:=W_{\Hilb}\times_{\Hilb^p}\Spec(\kappa(x))\in\Hilbf^p_K(\Spec(\kappa(x)))$. Furthermore let $\sheafi^{(x)}\subset\sheafo_{\proj[\kappa(x)]}$ denote the associated coherent ideal sheaf. It has Hilbert polynomial $q(t):=\binom{t+n}{n}-\nolinebreak p(t)$.

\begin{definition}
For $i\in\N$ and $x\in\Hilb^p$ define the $i$th \emph{cohomological Hilbert function}
$$h^i_{x}:\Z\ra\N,\ m\mapsto\dim_{\kappa(x)}H^i\bigl(\proj[\kappa(x)],\sheafi^{(x)}(m)\bigr).$$
\end{definition}

Fix a sequence $(f_i)_{i\in\N}$ of numerical functions $f_i:\Z\ra\N$. Then by the Semicontinuity Theorem 
$$H^{\geq}:=H^{\geq}_{\proj}:=\{x\in\Hilb^p\mid h^i_x\geq f_i\ \forall\,i\in\N\}$$ 
is a closed subspace of $\Hilb^p$ and 
$$H^{=}:=H^{=}_{\proj}:=\{x\in\Hilb^p\mid h^0_x=f_0,\ h^i_x\geq f_i\ \forall\,i\geq 1\}$$ 
is a locally closed subspace of $\Hilb^p$. 

In the rest of this section we prove that $H^{\geq}$ and $H^{=}$ are connected. Therefore, we redefine the polynomial ring $S$ by $S:=K[X_0,\dots,X_n]$, where $K$ is a field of characteristic zero as before.

\begin{notation}
Define
$$\Ideals^=:=\Big\{\ideala\subset S\Big|\ \parbox{6cm}{$\ideala$ is a saturated homogeneous ideal with $h_\ideala=f_0$ and $h^i_\ideala\geq f_{i-1}$ for all $i\geq 2$}\Big\},$$
$$\Ideals^\geq:=\Big\{\ideala\subset S\Big|\ \parbox{8.3cm}{$\ideala$ is a saturated homogeneous ideal with Hilbert polynomial $q$ such that $h_\ideala\geq f_0$ and $h^i_\ideala\geq f_{i-1}$ for all $i\geq 2$}\Big\},$$
where $h^i_\ideala:\Z\ra\N,\ m\mapsto\dim_K H^i_{S_+}(\ideala)_m$ denotes the $i$-th locally cohomological Hilbert function of $\ideala$.
\end{notation}

\begin{proposition}\label{connectedness by Groebner def}
The sets $\,\Ideals^=$ and $\,\Ideals^\geq$ are connected by Gr\"obner deformations.
\end{proposition}

\begin{proof}
Let $\ideala\in\Ideals^=$ and apply Proposition \ref{Mall's main theorem} to find two sequences of Mall ideals $\idealc_1$,~\dots,~$\idealc_r$ and $\ideald_1$, \dots, $\ideald_s$ with the appropriate properties. 

By \cite{S} we know that $h^i_\idealc\leq h^i_{(\init_\deglex\idealc)^\sat}$ for all $i\geq 2$ and all homogeneous ideals $\idealc\subset S$. Since Mall ideals are sequentially Cohen-Macaulay (Theorem \ref{Mall ideals are sCM}), we know by \cite{HS} and (Theorem \ref{Gin = in for Mall ideals}) that $h^i_\idealc=h^i_{(\init_\degrevlex\idealc)^\sat}$ for all $i\geq 2$ and all Mall ideals $\idealc\subset S$. Hence we get the situation

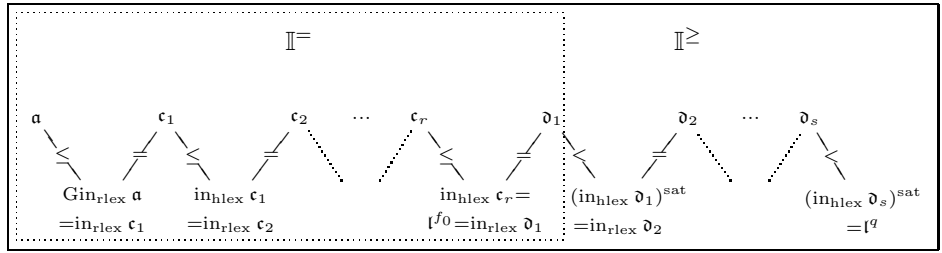
\begin{center}
\framebox{$\def\objectstyle{\scriptstyle}\def\labelstyle{\scriptstyle}\xymatrix@C=-2pt@R=17pt{
*i{\txt{100\\2}}&&&&\txt{\ \\$\Ideals^=$}&&&&&&\txt{\ \\$\Ideals^\geq$}\\
\ideala\ar@{-}[dr]|-{\leq}
&&\idealc_1\ar@{-}[dl]|-{=}\ar@{-}[dr]|-{\leq}
&&\idealc_2\ar@{-}[dl]|-{=}\ar@{.}[dr]
&\dots
&\idealc_r\ar@{.}[dl]\ar@{-}[dr]|-{\leq}
&&\ideald_1\ar@{-}[dl]|-{=}\ar@{-}[dr]|-{<}
&&\ideald_2\ar@{-}[dl]|-{=}\ar@{.}[dr]
&\dots
&\ideald_s\ar@{.}[dl]\ar@{-}[dr]|-{<}\\
&\txt{$\scriptstyle\Gin_\degrevlex\ideala$\\$\scriptstyle=\init_\degrevlex\idealc_1$}
&&\txt{$\scriptstyle\init_\deglex\idealc_1$\\$\scriptstyle=\init_\degrevlex\idealc_2$}
&&*+<8ex,5.2ex>{}
&&\txt{$\scriptstyle\init_\deglex\idealc_r=$\\\hspace*{-7pt}$\scriptstyle\ideall^{f_0}=\init_\degrevlex\ideald_1$\hspace*{-7pt}}
&*+<2ex,5.2ex>{}
&\txt{\hspace*{-6pt}$\scriptstyle(\init_\deglex\ideald_1)^\sat$\hspace*{-14pt}\\$\scriptstyle=\init_\degrevlex\ideald_2$}
&&*+<8ex,5.2ex>{}
&&\txt{\hspace*{-10pt}$\scriptstyle(\init_\deglex\ideald_s)^\sat$\hspace*{-10pt}\\$\scriptstyle=\ideall^q$}
\save "1,1"."3,9"*[F.]\frm{}\restore
}$\hspace{10pt}}
\end{center}
where the signs indicate the following:
\sk
\begin{list}{}{\setlength{\labelwidth}{6.0ex} \setlength{\leftmargin}{7.5ex} \setlength{\labelsep}{1.5ex}
\setlength{\itemindent}{0ex} \setlength{\listparindent}{0ex}
\setlength{\parsep}{0ex} \setlength{\itemsep}{0ex}
\setlength{\topsep}{0ex} \setlength{\partopsep}{0ex}
\setlength{\parskip}{0ex}}
\item[``$=$'':] The Hilbert function and the locally cohomological Hilbert functions remain constant. 
\item[``$\leq$'':] The Hilbert function remains constant and the locally cohomological Hilbert functions do not decrease.
\item[``$<$'':] The Hilbert function increases and the locally cohomological Hilbert functions do not decrease.
\end{list}
\end{proof}

\begin{theorem}\label{connectedness by lines}
The sets $H^{\geq}$ and $H^{=}$ are connected by lines.
\end{theorem}

\begin{proof}
Assume at first that the field $K$ is algebraically closed. Then the closed points of $\Hilb^p$ are precisely the saturated homogeneous ideals of $S$ with Hilbert polynomial $q$. By the Serre-Grothendieck Correspondence the set of closed points of $H^{=}$ equals $\Ideals^=$ and the set of closed points of $H^{\geq}$ equals $\Ideals^\geq$. By the previous Proposition, the only fact to prove is the following:

\sk
\emph{Claim: Let $\ideala\subset S$ be a homogeneous ideal and $\tau$ a term order of $S$ such that $\ideala$, $\init_\tau\ideala\in\Ideals_=$ or such that $\ideala$, $\init_\tau\ideala\in\Ideals_\geq$. Then $\ideala$ and $\init_\tau\ideala$ are connected by a line which lies entirely in $\Ideals_=$, $\Ideals_\geq$ respectively.}
\sk

In fact, by means of weight orders, described for example in Chapter 15.8 of \cite{E}, one constructs an ideal $\beta^z(\ideala)\subset S[z]$ and a flat family $S[z]/\beta^z(\ideala)$ over $K[z]$ whose fibre over $0$ is $S/\init_\tau\ideala$, and whose fibre over $(z-u)$ for $u\in K\setminus\{0\}$ is isomorphic to $S/\ideala$. A detailed scheme theoretic proof of the Claim may be found in section 3.3 of \cite{F}.

Now, assume that $K$ is an arbitrary field of characterisic zero. Endow the locally closed subspaces $H^{\geq}_{\proj}$, $H^{=}_{\proj}\subset\Hilb^p_{\proj}$ with the induced reduced scheme structure. Let $k$ be the algebraic closure of $K$. The proof is finished if we show:

\sk
\emph{Claim: $H^{\geq}_{\proj[k]}\cong(H^{\geq}_{\proj}\times_{\Spec(K)}\Spec(k))_\red$ and $H^{=}_{\proj[k]}\cong(H^{=}_{\proj}\times_{\Spec(K)}\Spec(k))_\red$.}

\sk
A detailed proof of this fact may be found in section 3.4 of \cite{F}. Its idea is the following: Define functors  $$\Hf^=_{\proj}:\sch_K^\red\ra\sets,\quad\Hf^\geq_{\proj}:\sch_K^\red\ra\sets$$ 
by assigning to each reduced $K$-scheme $T$ the sets 
$$\Hf^=_{\proj}(T):=\left\{W\in\Hilbf^p_{\proj}(T)\Bigg|\ \parbox{7.2cm}{\small{If $g:T\ra\Hilb^p$ is a morphism such that $W$ is the pull back of the universal subscheme $H_{\Hilb}$ by $g$, then $h^0_{g(x)}=f_0$ and $h^i_{g(x)}\geq f_i\ \forall i\geq 1\ \forall x\in T$.}}\right\},$$
$$\Hf^\geq_{\proj}(T):=\left\{W\in\Hilbf^p_{\proj}(T)\Bigg|\ \parbox{7.2cm}{\small{If $g:T\ra\Hilb^p$ is a morphism such that $W$ is the pull back of the universal subscheme $H_{\Hilb}$ by $g$, then $h^i_{g(x)}\geq f_i\ \forall i\in\N\ \forall x\in T$.}}\right\},$$
respectively. The functors $\Hf^=_{\proj}$, $\Hf^\geq_{\proj}$ are subfunctors of the Hilbert functor. It is easily shown that  $\Hf^=_{\proj}$, $\Hf^\geq_{\proj}$ are represented by the schemes $H^{\geq}_{\proj}$, $H^{=}_{\proj}$, respectively. Now, the claim follows from general nonsense  (cf. \cite[0.1.3.10]{EGA I}).
\end{proof}
\section*{Acknowledgements}
I wish to thank my thesis advisor Prof.~M.~Brodmann for his continuous support, E.~Sbarra for stimulating discussions, and Prof.~B.~Sturmfels for helpful comments.


\begin{thebibliography}{99}
\bibitem[1]{AS} K. Altmann and B. Sturmfels: The Graph of Monomial Ideals, \emph{Journal of Pure and Applied Algebra} 201 (2005), 250--263.

\bibitem[2]{BS} M. Brodmann and R.Y. Sharp: \emph{Local cohomology -- an algebraic introduction with geometric applications,} Cambridge studies in advanced mathematics 60, Cambridge University Press (1998).

\bibitem[3]{BH} W. Bruns and J. Herzog: \emph{Cohen Macaulay rings,} Cambridge studies in advanced mathematics 39, Cambridge University Press (1993).

\bibitem[4]{Cocoa}CoCoATeam: CoCoA: a system for doing Computations in Commutative Algebra, available at http://cocoa.dima.unige.it.

\bibitem[5]{C} A. Conca: Koszul homology and extremal properties of Gin and Lex, \emph{Trans.\ Amer.\ Math.\ Soc.}\ 356 (2004), 2945--2961.

\bibitem[6]{E} D. Eisenbud: \emph{Commutative algebra with a view towards algebraic geometry,} Graduate Texts in Mathematics, Springer-Verlag New York (1995).

\bibitem[7]{F} S. Fumasoli: Connectedness of Hilbert scheme strata defined by bounding cohomology, Ph.\,D. Thesis, Universit\"at Z\"urich (2005). arXiv:math.AC/0509123.

\bibitem[8]{Go1988} G. Gotzmann: Durch Hilbertfunktionen definierte Unterschemata des Hilbertschemas, \emph{Comment.\ Math.\ Helvetici} 63 (1988), 114--149.

\bibitem[9]{Gr} A. Grothendieck: Techniques de construction et th\'eor\`emes d'exis\-tence en g\'eom\'etrie alg\'ebri\-que IV. Les sch\'emas de Hilbert, \emph{S\'emi\-naire Bourbaki} 221 (1961).

\bibitem[10]{EGA I} A. Grothendieck and J. Dieudonn\'e: \emph{\'El\'ements de G\'eom\'etrie Alg\'ebri\-que~I,} Die Grund\-lehren der mathematischen Wissenschaf\-ten, Sprin\-ger-Ver\-lag Berlin Heidelberg (1971).

\bibitem[11]{H1966} R. Hartshorne: Connectedness of the Hilbert Scheme, \emph{Publ.\ Math.\ IHES} 29 (1966), 261--304.

\bibitem[12]{HPV} J. Herzog, D. Popescu, and M. Vladoiu: On the Ext-modules of ideals of Borel type, in: L. L. Avramov, M. Chardin, M. Morales, and C. Polini (Eds.), \emph{Commutative Algebra. Interactions with Algebraic Geometry,} Contemp.\ Math.\ 331, Amer.\ Math.\ Soc., Providence, RI, (2003), 171--186.

\bibitem[13]{HS} J. Herzog and E. Sbarra: Sequentially Cohen-Macaulay modules and local Cohomology, in: R. Parimala (Ed.), \emph{Algebra, arithmetic and geometry, Part I, II (Mumbai, 2000),} Tata Inst.\ Fund.\ Res., Bombay (2002), 327--340.

\bibitem[14]{M1997} D. Mall: Betti numbers, Castelnuovo Mumford regularity, and generalisations of Macaulay's Theorem, \emph{Communications in Algebra} 25(12), (1997), 3841--3852.

\bibitem[15]{M2000} D. Mall: Connectedness of Hilbert function strata and other connectedness results, \emph{Journal of Pure and Applied Algebra} 150 (2000), 175--205.

\bibitem[16]{P} K. Pardue: Deformations of graded modules and connected loci on the Hilbert scheme, The Curves Seminar at Queen's, Vol.\ XI, \emph{Queen's Papers in Pure and Appl.\ Math.}\ 105 (1997), 132--149.

\bibitem[17]{PS} I. Peeva and M. Stillman: Connectedness of Hilbert schemes, \emph{J. Algebraic Geom.}\ 14  (2005), 193--211.

\bibitem[18]{S} E. Sbarra: Upper bounds for local cohomology for rings with given Hilbert function, \emph{Communications in Algebra} 29(12), (2001), 5383--5409.

\bibitem[19]{Sp} E. Sperner: \"Uber einen kombinatorischen Satz von Macaulay und seine Anwendungen auf die Theorie der Polynomideale, \emph{Abh.\ math.\ Sem.\ Univ.\ Hamburg} 7 (1930), 149--163.

\bibitem[20]{St} S. A. Str{\o}mme:  Elementary introduction to representable functors and Hilbert schemes,  in P. Pragacz (Ed.), \emph{Parameter Spaces,} Banach Center Publications 36, Warszawa (1996), 179--198.

\bibitem[21]{Stu} B. Sturmfels: \emph{Gr\"obner bases and convex polytopes}, AMS University Lecture Series 8, Amer.\ Math.\ Soc., Providence, RI, (1996).
\end{thebibliography}
\end{document}